\documentclass[10pt]{amsart} 
\usepackage{amssymb}

\theoremstyle{plain}

\newtheorem{proposition}{Proposition}[section]
\newtheorem{corollary}{Corollary}[section]
\newtheorem{lemma}{Lemma}[section]

\theoremstyle{definition}
\newtheorem{definition}{Definition}[section]

\newtheorem{remark}{Remark}[section]

\newcommand{\sn}{\mathbb S}

\DeclareMathOperator{\ricci}{Ricci}
\DeclareMathOperator{\Id}{Id}
\DeclareMathOperator{\tr}{trace}

\begin{document}
\title{The harmonicity of nearly cosymplectic structures}
\author{E. Loubeau}
\address{D{\'e}partement de Math{\'e}matiques \\
Universit{\'e} de Bretagne Occidentale \\
6, avenue Victor Le Gorgeu \\
CS 93837, 29238 Brest Cedex 3, France}
\email{Eric.Loubeau@univ-brest.fr}

\author{E. Vergara-Diaz}
\address{School of Mathematics, Trinity College Dublin, Dublin 2, Ireland}
\email{evd@maths.tcd.ie}

\thanks{Research carried out under the EC Marie Curie Action no. 219258}

\begin{abstract}
Almost contact structures can be identified with sections of a twistor bundle and this allows to define their harmonicity, as sections or maps. We consider the class of nearly cosymplectic almost contact structures on a Riemannian manifold and prove curvature identities which imply the harmonicity of their parametrizing section, thus complementing earlier results on nearly-K{\"a}hler almost complex structures.
\end{abstract}

\maketitle

\section{Introduction}

The concept of harmonic maps is probably best seen as a variational problem on the infinite dimensional manifold of maps between two Riemannian manifolds and, as such, provides a (natural) criterion, or sieve, to select remarkable mappings and the recurring interactions between geometry, or even topology, and analysis of the early results reinforce this impression. 

Paradoxically, vector fields are omnipresent in Differential Geometry but seldom considered as maps, most likely because there is no canonical metric on the tangent bundle and constructing one opens up a wide debate on its traits. If one opts for the Sasaki metric, clearly the first at hand, the harmonic map problem and the weaker harmonic vector fields, obtained by only admitting variations through sections, prove to be disappointing as a compact domain will force harmonicity to be exactly parallelism. When possible, one can further ease up these restrictions working with the unit tangent bundle and this turns out to be more successful as, for example, the Hopf vector field on $\sn^3$ is a unit harmonic vector field (and map to boot), cf. the survey~\cite{M}. 

By the same token, this theoretical framework can be put into practice for any fibre bundle yielding the notions of harmonic sections and harmonic maps. This was developed for the twistor space (and its Sasaki-type metric) of an almost Hermitian manifold in~\cite{CMW1,CMW2}, and via a relatively sophisticated construction of projections, the equation characterizing harmonic sections was given in the more palatable language of tensor fields as
$$[\nabla^{*}\nabla J , J]=0 ,$$
with $\nabla^{*}\nabla = - \tr \nabla^{2}$, while the more restrictive harmonic map equation requires the extra condition
$$ \langle \nabla_{E_i} J , [R(E_{i}, X),J]\rangle =0,$$
for any tangent vector $X$ and an orthonormal frame $\{E_{i}\}_{i=1,\dots,2n+1}$.

The analogous notion on odd-dimensional manifolds is the almost contact structures, a slightly more demanding geometric object since essentially constituted of a tensor and a vector field but often defined on a Riemannian manifold $(M,g)$ as a $(1,1)$-tensor $\theta$, a 1-form $\eta$ and a unit vector field $\xi$ such that
$$\theta^2 = -\Id + \eta\otimes \xi ; \quad \eta(\xi)=1.$$
The choice of the metric $g$ being paramount to all our considerations, we take the blanket assumption on the metric compatibility of almost contact structures.
The existence of such a triple is in fact a reduction to $U(n)\times U(1)$ of the structure group and, from this viewpoint, an analysis similar to the case of almost complex structures is carried out in~\cite{VW1}, where the \emph{first} and \emph{second} harmonic section equations are computed to be
\begin{equation}\label{hse1}
[\bar{\nabla}^{*}\bar{\nabla} J , J]=0 ,
\end{equation}
and
\begin{equation}\label{hse2}
\nabla^{*}\nabla \xi = |\nabla\xi|^2 - (1/2) J\circ \tr(\bar{\nabla}J\otimes\nabla\xi),
\end{equation}
and the addition of 
\begin{equation}\label{hme}
\langle \bar{\nabla}_{E_i} J , [\bar{R}(E_{i}, X),J]\rangle + 8\langle \nabla_{E_i}\xi , R(E_{i},X)\xi \rangle, \quad \forall X\in TM,
\end{equation}
makes up the whole tension field (cf. the end of the section for notations).

The various classes of almost contact structures can then be studied under the light of harmonicity and, for instance, a nearly cosymplectic manifold with parallel characteristic vector field or the standard cosymplectic structure on the $5$-sphere defines a harmonic map~\cite{VW1,VW2}. A comparison should be drawn with its even-dimensional counterpart since nearly-K{\"a}hler almost complex structures are always harmonic~\cite{CMW1}, due to a curvature property of Gray~\cite{G2}. 

The aim of this article is to prove that a nearly cosymplectic almost contact structure always defines a harmonic map, hence a harmonic section. This relies on a series of curvature identities obtained from the second covariant derivative of the tensor field $\theta$ and the symmetries of $\nabla\theta$ and the Riemannian curvature tensor, but not before a re-writing of the harmonic section equations in a more amenable form, the harmonic map equation being from the start curvature compatible. Once this is achieved the conclusion follows more or less straightforwardly. This approach owes a lot to~\cite{G1} and~\cite{CMW2}.

For an almost contact structure $(\theta, \eta,\xi)$ on the Riemannian manifold $(M,g)$, we denote by $\mathcal F$ the distinguished distribution, called \emph{horizontal}, complementary to the $\xi$-direction, by $J$ the restriction of $\theta$ to $\mathcal F$ and $\bar{\nabla}$ the projection onto $\mathcal F$ of the Levi-Civita connection of $(M,g)$. The set of vector fields $\{F_i\}_{i=1,\dots,2n}$ will be a local orthonormal frame of the distribution $\mathcal F$, sometimes replaced, for computational reasons, by $\{\theta F_i\}_{i=1,\dots,2n}$, and completed by $\xi$ into the orthonormal frame $\{E_i\}_{i=1,\dots,2n+1}= \{F_i\}_{i=1,\dots,2n}\cup \{\xi\}$ of $TM$. One tensor, in particular, spontaneously appears in our computations as a bridge between harmonicity and curvature in presence of an almost complex structure.

\begin{definition}~\cite{G2} Let $(M,g)$ be a Riemannian manifold with curvature tensor $R$. If $(\theta, \xi,\eta)$ is an almost contact structure, we define the Ricci-$*$ curvature operator by
$$\ricci^* (X,Y)= \sum_{i=1}^{2n} g(R(X,F_i)\theta F_i, \theta Y),$$
where $\{F_i \}_{i=1,\dots,2n}$ is an orthonormal basis of the $\mathcal F$-distribution.
\end{definition}

We will adopt the Einstein convention on the summation of repeated indices and the following sign for the Riemannian curvature tensor:
$$R(X,Y)= [\nabla_{X},\nabla_{Y}] - \nabla_{[X,Y]}Z ,$$
with
$$R(X,Y,Z,W)= g(R(X,Y)Z,W).$$

\section{Harmonic section equations and Ricci-$*$}

This part lays the ground for the main computations of the article, by transforming the harmonic section equations into curvature conditions. The relevant operator for the first set of equations is a Ricci-type curvature introduced by Gray, since second-order covariant derivatives of $J$ will involve the $\mathcal F$-component of the Riemann curvature tensor and when we trace its Lie bracket with $J$, Ricci-$*$ naturally appears. Meanwhile, exploiting the symmetries of $\theta$, the second harmonic section equation is re-written as a curvature expression for the rough Laplacian of the Killing vector field $\xi$.

\begin{definition}
A \emph{nearly cosymplectic structure} is an almost contact structure $(M^{2n+1},\theta, \xi,\eta,g)$ such that
$$(\nabla_X\theta) (Y)+ (\nabla_Y \theta) (X)=0 ,$$
and we then say that $(M^{2n+1},\theta, \xi,\eta,g)$ is a {\em nearly cosymplectic manifold}.
\end{definition}

This clearly implies that the unit vector field is Killing and has geodesic integral curves.

We now make the blanket assumption that $(M^{2n+1},\theta, \xi,\eta,g)$ is a nearly cosymplectic manifold and the next lemma collects some of their more immediate properties.

\begin{lemma} \label{nc1}\cite{Blair}
For any vectors $X$ and $Y$ in the $\mathcal{F}$-distribution, we have
\begin{align*}
&(\bar{\nabla}_X J)(Y)+ (\bar{\nabla}_Y J)(X)=0 ;\\
&(\bar{\nabla}_{\xi}J) (X)=J(\nabla_X \xi) = \theta (\nabla_X \xi ) = (\nabla_{\xi} \theta) (X) ;\\
&J(\bar{\nabla}_{\xi}J)(X)=-\nabla_X \xi=\theta (\nabla_{\xi}\theta) (X) ;\\
&(\bar{\nabla}_X J)(Y)= -(\bar{\nabla}_{JX }J)(JY) .
\end{align*}
\end{lemma}
Combining this last equality and the fact that $J$ preserves orthonormal bases of $\mathcal{F}$, shows that $J$ is divergence free, i.e. $\bar{\delta}J=(\bar{\nabla}_{F_i}J)(F_i)$ vanishes, $\{F_i\}_{i=1,\dots,2n}$ being an orthonormal basis of the horizontal distribution $\mathcal{F}$.

The covariant derivation of the defining equation of $\theta$ leads to a generalization of the last equation of Lemma~\ref{nc1}.
\begin{lemma}\label{nc7b}
Let $X$ and $Y$ be vectors tangent to $M$. Then 
$$(\nabla_{\theta X} \theta) (\theta Y) = -(\nabla_X \theta) (Y) +g(Y,\xi)\nabla_{\theta X} \xi + g(X,\xi)\theta (\nabla_Y \xi).$$
\end{lemma}
\begin{proof} 
Let $X$ and $Y$ be tangent vectors, extended locally into vector fields, then
\begin{align*}
(\nabla_{\theta X} \theta) (\theta Y)=& \nabla_{\theta X}(\theta^2 Y)  - \theta (\nabla_{\theta X} \theta) (Y)-\theta^2 (\nabla_{\theta X} Y)\\
= & g(Y, \nabla_{\theta X} \xi)\xi + \eta(Y)\nabla_{\theta X} \xi -\theta (\nabla_{\theta X} \theta) (Y)\\
= & g(Y, \nabla_{\theta X} \xi)\xi + \eta(Y)\nabla_{\theta X} \xi +  \theta (\nabla_{Y} \theta ) (\theta X)\\
= & g(Y, \nabla_{\theta X} \xi)\xi + \eta(Y)\nabla_{\theta X} \xi -\theta (\nabla_Y X)\\
& \theta (g(\nabla_Y X,\xi)\xi + g(X,\nabla_Y \xi )\xi) +\eta(X)\theta(\nabla_Y \xi)\\
&+ (\nabla_Y \theta) (X) +\theta (\nabla_Y X)+g( \theta X, \nabla_Y\xi )\xi\\
= & \eta(Y)\nabla_{\theta X} \xi + \eta(X) \theta (\nabla_Y \xi)-(\nabla_X\theta) (Y),
\end{align*}
since $\xi$ is a Killing vector field.
\end{proof}

\begin{remark}\label{nc8} 
Setting $Y=\xi$ in Lemma \ref{nc7b} we obtain, for any $X\in TM$,
$$ \theta (\nabla_X \xi) = -\nabla_{\theta X} \xi ,$$
since $\xi$ has geodesic integral curves.
\end{remark}
One can also obtain the commutation of $\theta$ and its derivative, from the relation $\theta^2 = -\Id +\eta\otimes\xi$.
\begin{lemma}\label{nca}
Let $X$ and $Y$ be vectors in $TM$. Then 
$$(\nabla_{X} \theta) (\theta Y) = -\theta(\nabla_X \theta) (Y) +g(Y,\nabla_X \xi) \xi + g(Y,\xi)\nabla_X\xi.$$
\end{lemma}
\begin{proof} 
Extend the vectors $X$ and $Y$ into local horizontal vector fields, then
\begin{align*}
(\nabla_{X} \theta) (\theta Y) &= \nabla_X (\theta^2 Y) - \theta (\nabla_X (\theta Y))\\
&=\nabla_X (-Y + \eta (Y)\xi) -\theta (\nabla_X \theta) (Y) -\theta^2 (\nabla_X Y)\\
&=-\theta (\nabla_X \theta) (Y) +g(Y,\nabla_X \xi) \xi + g(Y,\xi)\nabla_X\xi .
\end{align*}
\end{proof}

\begin{remark}
If, in the above equation, one chooses $X=Y$ in $\mathcal{F}$, then $(\nabla_X \theta)(\theta X)$ must vanish, since $\xi$ is Killing and $\theta$ nearly cosymplectic.
\end{remark} 
In order to compute the harmonic section equations, we need the commutator of the restriction to $\mathcal F$ of $\theta$ and its rough Laplacian.
\begin{lemma}\label{nc6}
Let $E$ be a vector in the distribution $\mathcal{F}$. Then
$$[\bar{\nabla}^2_{E,E}J,J]=-[\bar{\nabla}^2_{JE,JE}J,J]+2[\bar{R}(E,JE),J],$$ 
where $\bar{R}$ is the curvature of $(\mathcal{F},\bar{\nabla})$.
\end{lemma}

\begin{proof} Let $x$ be a point in $M$ and $E$ and $X$ vectors in $\mathcal{F}_x$ which we extend to local sections of $\mathcal{F}$
such that $\bar{\nabla}X=\bar{\nabla}E=0$ at the point $x$.
First, note that, using Lemma~\ref{nc1},
\begin{align} 
[E, \theta E]&=[E,JE]\label{clubsuit}\\
& = \bar{\nabla}_E (JE)-g( \nabla_E (JE),\xi ) \xi+g( \nabla_{JE} E,\xi )\xi \notag\\
& = (\bar{\nabla}_{E}J)(E)+g( JE, \nabla_E \xi ) \xi - g( E, \nabla_{JE} \xi ) \xi \notag\\
& = g( JE, \nabla_E \xi ) \xi + g( E, J(\bar{\nabla}_{\xi} J)(JE) ) \xi \notag\\ 
& = g( JE, \nabla_E \xi ) \xi + g( E, (\bar{\nabla}_{\xi} J)(E) ) \xi \notag\\
& = - g( E, J(\nabla_E \xi) ) \xi + g( E,J (\nabla_E \xi) ) \xi\notag\\
&=0 .\notag
\end{align}
Now, by the Leibniz rule
$$(\bar{\nabla}_E\bar{\nabla}_EJ)(JX)=\bar{\nabla}_E(\bar{\nabla}_EJ (JX))-(\bar{\nabla}_{E}J)\circ (\bar{\nabla}_E J) (X),$$
and the first term may be expanded as follows
\begin{align*}
\bar{\nabla}_E((\bar{\nabla}_EJ) (JX))& =  \bar{\nabla}_E((\bar{\nabla}_{JE}J) (X))\\
& =\bar{\nabla}_E\bar{\nabla}_{JE}(JX)-\bar{\nabla}_E(J(\bar{\nabla}_{JE}X))\\
& =\bar{\nabla}_E\bar{\nabla}_{JE}(JX)-J(\bar{\nabla}_E\bar{\nabla}_{JE}X)\\
& =\bar{\nabla}_{JE}\bar{\nabla}_E(JX)+ \bar{\nabla}_{[E,JE]}(JX)+\bar{R}(E,JE)(JX)\\
& -J\bar{\nabla}_{JE}\bar{\nabla}_{E}X-J\bar{\nabla}_{[E,JE]}X-J\bar{R}(E,JE)X\\
& =\bar{\nabla}_{JE}\bar{\nabla}_E(JX)-J\bar{\nabla}_{JE}\bar{\nabla}_E X+[\bar{R}(E,JE),J]X .
\end{align*}
Therefore
\begin{align*}
&\bar{\nabla}_E((\bar{\nabla}_E J)(JX))-[\bar{R}(E,JE),J]X\\ &=\bar{\nabla}_{JE}((\bar{\nabla}_EJ)(X))+\bar{\nabla}_{JE}(J(\bar{\nabla}_EX))-J\bar{\nabla}_{JE}\bar{\nabla}_EX\\
&=-\bar{\nabla}_{JE}((\bar{\nabla}_{JE}J)(JX))\\
&=-(\bar{\nabla}_{JE}\bar{\nabla}_{JE}J)(JX)-(\bar{\nabla}_{JE}J)(\bar{\nabla}_{JE}(JX))\\
&=-(\bar{\nabla}_{JE}\bar{\nabla}_{JE}J)(JX)-(\bar{\nabla}_{JE}J)\circ(\bar{\nabla}_{JE}J)(X),
\end{align*}
and it follows that
\begin{align}\label{eq1}
(\bar{\nabla}_E\bar{\nabla}_EJ) (JX)& = 
-(\bar{\nabla}_{JE}\bar{\nabla}_{JE}J)(JX)-(\bar{\nabla}_{JE}J)\circ (\bar{\nabla}_{JE}J)(X)\\
&-(\bar{\nabla}_{E}J)\circ (\bar{\nabla}_EJ)(X) +[\bar{R}(E,JE),J]X \notag.
\end{align}
Again, by Lemma~\ref{nc1}
\begin{align*}
(\bar{\nabla}_{JE}J)\circ (\bar{\nabla}_{JE}J)(X)&=(\bar{\nabla}_EJ)\circ J\circ (\bar{\nabla}_EJ)(JX)\\
&=(\bar{\nabla}_EJ)\circ(\bar{\nabla}_EJ)(X) ,
\end{align*}
so Equation~\eqref{eq1} can be rewritten
\begin{align*}
(\bar{\nabla}_E\bar{\nabla}_EJ)(JX)& = -(\bar{\nabla}_{JE}\bar{\nabla}_{JE}J)(JX)-2(\bar{\nabla}_EJ)\circ (\bar{\nabla}_EJ)(X) +[\bar{R}(E,JE),J]X ,
\end{align*}
and since 
\begin{align*} 
(\bar{\nabla}^2_{E,E} J) (JX) &=( \bar{\nabla}_E \bar{\nabla}_E J) (JX) - (\bar{\nabla}_{\nabla_E E} J)(JX)\\
&= (\bar{\nabla}_E \bar{\nabla}_E J) (JX) - (\bar{\nabla}_{g( \nabla_E E ,\xi ) \xi} J)(JX)\\
&= (\bar{\nabla}_E \bar{\nabla}_E J) (JX) + (\bar{\nabla}_{g(  E ,\nabla_E \xi ) \xi} J)(JX)\\
&= (\bar{\nabla}_E \bar{\nabla}_E J)(JX), 
\end{align*}
it follows that
\begin{align*}
J(\bar{\nabla}_{E,E}^2J)(X)& = -J(\bar{\nabla}_{JE}\bar{\nabla}_{JE}J)(X)
-2(\bar{\nabla}_EJ)\circ (\bar{\nabla}_EJ)(X) -[\bar{R}(E,JE),J]X .
\end{align*}
In conclusion
\begin{align*}
[\bar{\nabla}^2_{E,E}J,J]& = -[\bar{\nabla}_{JE}\bar{\nabla}_{JE}J ,J]+2[\bar{R}(E,JE),J].
\end{align*}
\end{proof}

To obtain the whole of the rough Laplacian of $J$, we need to determine the contribution of $\xi$.
\begin{lemma} 
$$[\bar{\nabla}_{\xi,\xi}^2J , J ]=[\nabla_{\xi,\xi}^2\theta , \theta ]= 0 .$$
\end{lemma}

\begin{proof} 
Let $x$ be a point in $M$ and $X$ a vector in $\mathcal{F}_{x}$, we extend $X$ to a local vector field, still called $X$, such that $\bar{\nabla}X =0$ at the point $x$. First we prove that
$$(\nabla_{\xi,\xi}^2 \theta) (X) = \theta R(\xi, X)\xi .$$
As the unit vector field $\xi$ has geodesic integral curves
$$\nabla_{\xi} X = \bar{\nabla}_{\xi} X + g( \nabla_{\xi} X,\xi ) = -g( X,\nabla_{\xi}\xi ) =0,$$
and with the Leibniz rule, we have
\begin{align*}
(\nabla^2_{\xi,\xi} \theta) (X)& =\nabla_{\xi}((\nabla_{\xi}\theta) (X)) -(\nabla_{\nabla_{\xi}\xi}\theta) (X) -(\nabla_{\xi} \theta) (\nabla_{\xi} X)\\
& =\nabla_{\xi}((\nabla_{\xi}\theta) (X)) \\
&=\nabla_{\xi}(\theta (\nabla_X \xi))\\
&=(\nabla_{\xi} \theta) (\nabla_X \xi)+\theta (\nabla_{\xi} \nabla_X \xi )\\
&=\theta (\nabla_{\nabla_X \xi} \xi)+\theta (\nabla_{\xi} \nabla_X \xi )-\theta (\nabla_X \nabla_{\xi}\xi )-\theta (\nabla_{\nabla_{\xi}X}\xi)\\
&=\theta R(\xi, X)\xi .
\end{align*}
Second, since $\xi$ has geodesic integral curves,
$$(\nabla_{\xi} \theta) (X) = (\bar{\nabla}_{\xi} J) (X)$$
and
$$(\nabla_{\xi,\xi}^2 \theta) (X) = (\bar{\nabla}_{\xi,\xi}^2 J) (X).$$
Then, as $\xi$ is Killing with geodesic integral curves and using Lemma~\ref{nc1}
\begin{align*}
(\nabla_{\xi,\xi}^{2}\theta) (X)&=\theta R(\xi,X)\xi\\
&=-\theta (\nabla^2_{X,\xi} \xi)\\
&= J(\nabla_{\nabla_X \xi}\xi)\\
&=-J^2 ((\bar{\nabla}_{\xi} J)(\nabla_X \xi))\\
&= (\bar{\nabla}_{\xi} J)(\nabla_X \xi)\\
&= -(\bar{\nabla}_{\xi} J)\circ J\circ(\bar{\nabla}_{\xi}J)(X)\\
&= J\circ (\bar{\nabla}_{\xi}J) \circ (\bar{\nabla}_{\xi}J)(X).
\end{align*}
Finally
\begin{align*} 
[\bar{\nabla}_{\xi,\xi}^{2}J, J ]X&=J\circ (\bar{\nabla}_{\xi}J) \circ (\bar{\nabla}_{\xi}J)(JX)- J^2\circ(\bar{\nabla}_{\xi}J)\circ (\bar{\nabla}_{\xi}J) (X)\\
&= 0.
\end{align*}
\end{proof}

The previous lemmas enable us to compute the first harmonic section equation~\eqref{hse1} in terms of the Riemann curvature tensor and reduce it to a $\theta$-invariance of Ricci-$*$, in perfect accordance with the case of almost complex structures~\cite{CMW2}. As we will see in the next section, this condition is actually automatically satisfied for nearly cosymplectic structures.

\begin{proposition}\label{ric*} 
Let $M^{2n+1}$ be a nearly cosymplectic manifold. Then the first harmonic equation is verified if and only if 
$$\ricci^*(\theta X, \theta Y)=   \ricci^*(X,Y) ,$$ 
for $X$ and $Y$ in $\mathcal{F}$.
\end{proposition}

\begin{proof} 
Let $x$ be a point of $M$ and $\{F_i \}_{i=1,\dots,2n}$ a local orthonormal frame of the $\mathcal F$-distribution such that, at the point $x$, $\bar{\nabla} F_i=0$.
Taking traces in Lemma~\ref{nc6}, we have
$$
-[\bar{\nabla}^*\bar{\nabla}J,J] - [\bar{\nabla}_{\xi,\xi}^2J,J] = [\bar{\nabla}^*\bar{\nabla}J,J] + [\bar{\nabla}_{\xi,\xi}^2J,J] +2[\bar{R}(F_i,JF_i),J],
$$
hence
$$-[\bar{\nabla}^*\bar{\nabla}J,J]  = [\bar{R}(F_i,JF_i),J],$$
by the previous lemma, and the first harmonic section equation is satisfied if and only if 
$$[\bar{R}(F_i,JF_i),J]=0.$$ 

Since the projection onto the $\mathcal F$-bundle of the curvature tensor is linked to the curvature tensor of $({\mathcal F},\bar\nabla )$ by the equation~\cite{thesis}
$$\bar{R}(X,Y) \beta = R^{\mathcal{F}}(X,Y) \beta+r(\nabla_X \xi,\nabla_Y \xi)\beta ,$$
where
$$r(\nabla_X \xi,\nabla_Y \xi)\beta= g(  \nabla_Y \xi , \beta ) \nabla_X  \xi -g(  \nabla_X \xi , \beta ) \nabla_Y  \xi ,$$
we can use the skew-symmetry of $\theta$ and Remark~\ref{nc8} to obtain
\begin{align}
r(\nabla_{F_i} \xi,\nabla_{\theta F_i} \xi) \theta \beta &= g(  \nabla_{\theta F_i} \xi , \theta \beta ) \nabla_{F_i} \xi -g(  \nabla_{F_i}\xi , \theta \beta ) \nabla_{\theta F_i}  \xi \notag\\
& =-g( \nabla_{ F_i} \xi  , \beta ) \nabla_{F_i} \xi -g(  \nabla_{\theta F_i}\xi , \beta ) \nabla_{\theta F_i}  \xi \notag\\
& =\theta r( \nabla_{F_i} \xi, \nabla_{\theta F_i} \xi)\beta,\label{eqA}
\end{align}
and the vanishing of the first harmonic section equation is then also equivalent to 
$$[R^{\mathcal{F}}(F_i,JF_i),J]=0.$$
Finally, let $Z$ and $W$ be in $\mathcal{F}$, by Bianchi's first identity
\begin{align*}
g( \bar{R}(F_i,JF_i) Z,W )&=g(R(F_i ,J F_i )Z,W ) + g(r(\nabla_{F_i}\xi,\nabla_{JF_i}\xi)Z,W),\\
& =-2 g( R(Z , F_i )JF_i,W )+ g( r(\nabla_{F_i}\xi,\nabla_{JF_i}\xi)Z,W) ,
\end{align*}
and, from Equation~\eqref{eqA}, we deduce
\begin{align*}
g([\bar{R}(F_i,JF_i),J](Z),W)& =g( [R^{\mathcal{F}}(F_i,JF_i) , J ](Z),W)\\
& =-2g( R^{\mathcal{F}}(JZ,F_i)JF_i-JR^{\mathcal{F}}(Z,F_i)JF_i , W ) \\
& = -2g( R(JZ,F_i)JF_i, W) +2g(JR(Z,F_i)JF_i,W )\\
& =-2 g(R(JZ,F_i)JF_i, W) - 2g(R(Z,F_i)JF_i,JW ) \\
& =2\ricci^*(\theta Z, \theta  W)-2 \ricci^*(Z,W),
\end{align*}
which proves the proposition.
\end{proof}

\begin{remark} 
When $\eta$ is closed, the second harmonic section equation is automatically satisfied since $\xi$ is then parallel.
Furthermore, in this case, the identity~\cite{Blair}
$$R(\theta X, \theta Y,\theta Z, \theta W ) =  R(X,Y,Z, W )$$
implies that the first harmonic section equation is also verified as
\begin{align*} 
\ricci^* (\theta Z, \theta W)&= - R(\theta Z,F_i,\theta F_i, W)\\
&=-R(Z,\theta F_i ,F_i,\theta W)\\
&=-R(Z,\theta^2 F_i ,\theta F_i,\theta W)\\
&=\ricci^* (Z,W).
\end{align*}
\end{remark}

The key ingredient to rewrite the second set of equations is Remark~\ref{nc8}, as it allows us to trace over $F_i$ and $\theta F_i$ and, swapping them around, introduce curvature terms. Again, this condition will turn out to be valid for any nearly cosymplectic structure.

\begin{proposition}\label{nc10}
Let $M^{2n+1}$ be a nearly cosymplectic manifold. Then the second harmonic section equation~\eqref{hse2} is equivalent to 
$$ \nabla^{*}\nabla \xi -|\nabla\xi|^2 \xi =- \tfrac{1}{2} [R(F_i, \theta F_i),\theta ] \xi ,$$ 
where $\{F_i \}_{i=1,\dots,2n}$ is a local orthonormal frame of the $\mathcal F$-distribution.
\end{proposition} 

\begin{proof}
As previously, let $x$ be a point of $M$ and $\{F_i \}_{i=1,\dots,2n}$ a local orthonormal frame of the $\mathcal F$-distribution such that, at the point $x$, $\bar{\nabla} F_i=0$. Then, by Remark~\ref{nc8} and Equation~\eqref{clubsuit}
\begin{align*} 
\theta \circ \tr( \nabla \theta \circ  \nabla \xi) &= \theta \circ (\nabla_{F_i} \theta ) (\theta (\nabla_{\theta F_i} \xi))\\
& = \theta \nabla_{F_i} (\theta ^2 (\nabla_{\theta F_i}\xi)) - \theta^2 (\nabla_{F_i} (\theta (\nabla_{\theta F_i} \xi))) \\
& = -\theta (\nabla_{F_i} \nabla_{\theta F_i} \xi )+\nabla_{F_i} (\theta (\nabla_{\theta F_i} \xi))-g(\nabla_{F_i} (\theta (\nabla_{\theta F_i} \xi)),\xi)\xi\\
& = -\theta [\nabla_{F_i} \nabla_{\theta F_i} \xi - \nabla_{\theta F_i} \nabla_{F_i} \xi - \nabla_{[F_i, \theta F_i]} \xi] -\theta (\nabla_{\theta F_i} \nabla_{F_i} \xi) \\
&- \theta (\nabla_{[F_i, \theta F_i]} \xi) +\nabla_{F_i} \nabla_{F_i} \xi-g(\nabla_{F_i} \nabla_{ F_i} \xi,\xi)\xi\\
& = -\theta R(F_i,\theta F_i) \xi-\theta (\nabla_{\theta F_i} \nabla_{F_i} \xi) +\nabla_{F_i} \nabla_{ F_i} \xi-g(\nabla_{F_i} \nabla_{ F_i} \xi,\xi)\xi\\
& = -\theta R(F_i,\theta F_i) \xi-\theta (\nabla_{\theta F_i} (\theta (\nabla_{\theta F_i} \xi) ))-(\nabla^*\nabla \xi)^{\mathcal{F}}\\
& = -\theta R(F_i,\theta F_i) \xi -\theta \circ \tr(\nabla \theta \circ \nabla \xi)-\theta^2( \nabla_{\theta F_i}\nabla_{\theta F_i} \xi)-(\nabla^*\nabla \xi)^{\mathcal{F}}\\
& =-\theta R(F_i,\theta F_i) \xi -\theta \circ \tr(\nabla \theta \circ \nabla \xi)+\nabla_{\theta F_i}\nabla_{\theta F_i} \xi\\
& - g(\nabla_{\theta F_i} \nabla_{\theta F_i} \xi,\xi ) \xi -(\nabla^*\nabla \xi)^{\mathcal{F}}\\
& =-\theta R(F_i,\theta F_i) \xi -\theta \circ \tr(\nabla \theta \circ \nabla \xi)+\nabla_{\theta F_i, \theta F_i}^2 \xi \\
& - g(\nabla_{\theta F_i , \theta F_i}^2  \xi,\xi ) \xi 
-(\nabla^*\nabla \xi)^{\mathcal{F}} \\
& = -\theta R(F_i,\theta F_i) \xi -\theta \circ \tr(\nabla \theta \circ \nabla \xi)-2 (\nabla^*\nabla \xi)^{\mathcal{F}}, 
\end{align*}
since 
\begin{align*} 
\nabla_{\theta F_i} (\theta F_i) &= (\nabla_{\theta F_i} \theta) (F_i) + \theta (\nabla_{\theta F_i} F_i) \\
&=-\nabla_{F_i}(\theta^2 F_i)+\theta (\nabla_{F_i} (\theta F_i))\\
&=\nabla_{F_i} F_i +\theta (\nabla_{F_i} \theta) (F_i) +\theta^2 (\nabla_{F_i} F_i)\\
&=g(\nabla_{F_i} F_i,\xi )\xi \\
&=0,
\end{align*} 
because $\xi$ is Killing, and this implies that
$$\nabla_{\theta F_i,\theta F_i}^2 \xi  = \nabla_{\theta F_i} \nabla_{\theta F_i} \xi .$$
\end{proof}

\section{The harmonic section equations}

The tools of this section are simply the symmetries of the Riemann curvature tensor and the tensors $\theta$ and $\nabla\theta$. Since, as is customary, second covariant derivatives of $\theta$ involve the curvature, its characteristic equation gives a first expression for sectional curvatures (Proposition~\ref{28}), which we extend to the full $(0,4)$-tensor by a polarization argument (Proposition~\ref{32}). Then adequate choices of vectors quickly yield that the first and second harmonic section equations vanish for nearly cosymplectic structures.

Consider the $2$-form $\Theta (X, Y) = g(X,\theta Y)$ and denote by $R$ the curvature tensor of the Riemannian manifold $(M,g)$, then
$$R(W,X) (\Theta) (Y,Z) = -\Theta (R(W,X)Y,Z)-\Theta (Y,R(W,X)Z) ,$$
\begin{equation}\label{ncb}
\nabla^2\Theta (W,X,Y,Z) =(\nabla^2_{W,X} \Theta)(Y,Z)= g(Y,(\nabla^2_{W,X} \theta)(Z)),
\end{equation}
and therefore
$$R(W,X) (\Theta) (Y,Z)= (\nabla^2_{W,X} \Theta) (Y,Z)-(\nabla^2_{X,W} \Theta)(Y,Z).$$

\begin{lemma}\label{20} 
Let $X$ and $Y$ be vectors tangent to $M$. Then
\begin{align*}
R(X,Y) \Theta (X,\theta Y) &= \nabla^2 \Theta (X,Y,X,\theta Y)\\
&= R(X,Y, X,Y) - R(X,Y, \theta X, \theta Y) - R(X,Y, X, \eta(Y) \xi).
\end{align*}
\end{lemma}
\begin{proof} 
First, observe that
$$\nabla^2 \Theta (W,X,Y,Z)=-\nabla^2 \Theta(W,Y,Z,X)=-\nabla^2 \Theta(W,Y,X,Z),$$
since, by Equation~\eqref{ncb},
\begin{align*} 
\nabla^2 \Theta (W,X,Y,Z)&=  g( Y ,(\nabla_{W,X}^2 \theta) (Z))\\
& = g(Y, \nabla_W((\nabla_X \theta) (Z)) )- g(Y,(\nabla_{\nabla_W X}  \theta) (Z)) -g(Y,(\nabla_{X}\theta)(\nabla_W Z))\\
& = -g(Y, \nabla_W((\nabla_Z\theta) (X)) )+ g(Y,(\nabla_{Z}\theta) (\nabla_W X))+g(Y,(\nabla_{\nabla_W Z}  \theta) (X))\\
& = -\nabla^2\Theta( W,Y,Z,X),
\end{align*}
while
\begin{align*} 
&\nabla^2 \Theta (W,X,Y,Z)= -g((\nabla_{W,X}^2 \theta) (Y), Z)\\
& =-g(\nabla_W((\nabla_X \theta) (Y)),Z )+ g((\nabla_{\nabla_W X}  \theta) (Y), Z) +g((\nabla_X \theta) (\nabla_W Y),Z)\\
& = g(\nabla_W ((\nabla_Y \theta) (X)),Z )- g((\nabla_{Y}  \theta) (\nabla_W X), Z) -g((\nabla_{\nabla_W Y} \theta) (X),Z)\\ 
& = g((\nabla^2_{W,Y} \theta) (X),Z)\\
&=-\nabla^2\Theta( W,Y,X,Z).
\end{align*}
In particular, for all $W,X,Z\in TM$,
$$\nabla^2\Theta (W,X,X,Z)=0.$$
Combining the previous calculations, we obtain
\begin{align*}
&(R(X,Y) \Theta) (X,\theta Y)=\nabla^2 \Theta (X,Y,X,\theta Y) -\nabla^2 \Theta (Y,X,X,\theta Y) \\
& =-\Theta (R(X,Y)X,\theta Y ) - \Theta (X,R(X,Y) \theta Y)\\
& = g(R(X,Y) X,Y) -g(R(X,Y X, \eta(Y) \xi )-g(R(X,Y)\theta X ,\theta Y),
\end{align*}
since $\theta^2=-I+\eta \otimes \xi$. 
\end{proof}

\begin{proposition}\label{28} 
Let $M^{2n+1}$ be a nearly cosymplectic manifold. Then for any $X$ and $Y$ in $TM$
$$|(\nabla_X \theta) (Y)|^2 +g^{2}(Y,\nabla_X \xi ) = -R(X,Y, X, Y ) + R(X,Y, \theta X, \theta Y).$$
\end{proposition}
\begin{proof}
Using the first Bianchi identity and the curvature expression for the second covariant derivative of a Killing vector field, we easily obtain
\begin{equation}\label{27}
g(R(X,Y) X,\eta (Y)\xi) = \eta(Y) g(\nabla_{X,X}^2\xi, Y).
\end{equation}
On the other hand, for all $X,Y\in TM$, by Equation~\eqref{ncb} and Lemma~\ref{nca}, we have
\begin{align} 
&\nabla^2\Theta (X,X ,Y , \theta Y) = g(Y, (\nabla^2_{X,X } \theta) (\theta Y))\label{ncd}\\
& = g(Y, \nabla_X ((\nabla_X \theta)(\theta Y))) -g(Y, (\nabla_X \theta) (\nabla_X (\theta Y)))-g(Y,(\nabla_{\nabla_X X} \theta)(\theta Y))\notag\\
&= - g(Y, \nabla_X [\theta (\nabla_X \theta) (Y)-  g(Y,\nabla_X \xi)\xi-\eta(Y)\nabla_X\xi]) -g(Y, (\nabla_X \theta)\circ(\nabla_X \theta)(Y))\notag\\
&- g(Y, (\nabla_X \theta)\circ \theta (\nabla_X Y)) +g(Y,\theta (\nabla_{\nabla_X X } \theta) (Y)) -2\eta(Y)g(Y,\nabla_{\nabla_X X} \xi)\notag\\
& = -g(Y, (\nabla_X \theta)\circ (\nabla_X \theta) (Y))-g (Y , \theta \nabla_X ((\nabla_X \theta) (Y))) +g(\nabla_X Y, \nabla_X \xi)\eta(Y)\notag\\ 
& +g(Y,\nabla_X\nabla_X \xi)\eta(Y)+g^{2}(Y,\nabla_X \xi) +\eta(\nabla_X Y)g(Y,\nabla_X \xi)+g^{2}(Y,\nabla_X \xi)\notag\\
&+\eta(Y)g(Y,\nabla_X \nabla_X \xi) +|(\nabla_X \theta)(Y)|^2+g(Y, \theta (\nabla_X \theta)(\nabla_X Y))-\eta(Y)g(\nabla_X Y ,\nabla_X \xi)\notag\\
& -g(Y,\nabla_X \xi )\eta(\nabla_X Y)+g(Y, \theta (\nabla_{\nabla_X X} \theta) (Y))-2g(Y, \nabla_{\nabla_X X} \xi)\eta(Y)\notag\\
& = -g(Y,\theta(\nabla^2_{X,X} \theta)(Y))+2|(\nabla_X \theta)(Y)|^2 + 2g^{2}(Y,\nabla_X \xi)+ 2\eta(Y)g(Y,\nabla_{X,X}^2 \xi).\notag
\end{align}
Since $\theta^2 =-I+\eta\otimes \xi$, we have
\begin{align*} 
(\nabla_X \theta^2) (Y) &=\nabla_X ( -Y+ \eta(Y)\xi) +\nabla_X Y -\eta (\nabla_X Y)\xi\\
&=\eta(\nabla_X Y)\xi+ g(Y,\nabla_X\xi)\xi+\eta(Y)\nabla_X\xi-\eta (\nabla_X Y)\xi\\
&=g(Y,\nabla_X \xi )\xi + \eta(Y) \nabla_X \xi,
\end{align*}
therefore
\begin{align*} 
(\nabla_{X,X} \theta^2) (Y) &=\nabla_X ((\nabla_X \theta^2) (Y)) -(\nabla_X \theta^2)(\nabla_X Y) -(\nabla_{\nabla_X X} \theta^2) (Y)\\
&=\nabla_X[g(Y,\nabla_X \xi )\xi +\eta(Y) \nabla_X \xi)]-g(\nabla_X Y,\nabla_X \xi )\xi \\
&-\eta(\nabla_X Y) \nabla_X \xi-g(Y,\nabla_{\nabla_X  X } \xi )\xi - \eta(Y) \nabla_{\nabla_X X} \xi \\
&=g(\nabla_X Y,\nabla_X \xi)\xi + g(Y, \nabla_X \nabla_X \xi)\xi+2g(Y,\nabla_X \xi) \nabla_X \xi \\
&+\eta(\nabla_X Y)\nabla_X \xi+\eta(Y)\nabla_X\nabla_X \xi -g(\nabla_X Y,\nabla_X \xi )\xi \\
&-\eta(\nabla_X Y) \nabla_X \xi -g(Y,\nabla_{\nabla_X  X } \xi )\xi - \eta(Y) \nabla_{\nabla_X X} \xi \\
&=g(Y,\nabla^2_{X,X }\xi )\xi + 2g(Y, \nabla_X \xi) \nabla_X \xi+ \eta(Y) \nabla_{X,X}^2 \xi .
\end{align*}
The covariant derivation of $\theta^2$ yields
$$(\nabla^2_{X,X} \theta^2) (Y) =(\nabla^2_{X,X} \theta) (\theta Y)+  \theta \circ (\nabla^2_{X,X} \theta)(Y)+ 2 (\nabla_X\theta) \circ (\nabla_X \theta) (Y),$$
and since $\theta$ and $\nabla^2_{X,X} \theta$ are skew-symmetric
$$g((\nabla^2_{X,X} \theta^2) (Y),Y )  = 2g (\theta \circ (\nabla^2_{X,X} \theta)(Y), Y) + 2 g((\nabla_X\theta)\circ(\nabla_X \theta)(Y),Y).$$
Hence 
\begin{equation}\label{25} 
g(\theta (\nabla_{X,X}^2 \theta) (Y) ,Y) = |(\nabla_X \theta) (Y)|^2 + g(Y,\nabla^2_{X,X} \xi)\eta(Y)+ g^{2}(Y,\nabla_X \xi) .
\end{equation}
Equation~\eqref{ncb} means that 
$$g((\nabla^2_{X,X} \theta) (\theta Y),Y ) = \nabla^2 \Theta (X,X , Y, \theta Y),$$
and the skew-symmetry of $\theta$ and its covariant derivatives imply that 
$$g((\nabla^2_{X,X} \theta) (\theta Y),Y ) =g(Y, \theta (\nabla^2_{X,X} \theta) (Y)),$$
hence, by Equations~\eqref{25} and \eqref{ncd},
$$ \nabla^2 \Theta (X,X,Y,\theta Y) = |(\nabla_X \theta) (Y)|^2 +\eta(Y) g(Y,\nabla^2_{X,X}\xi) +g^{2}(Y,\nabla_X \xi) .$$
Moreover, we have  
$$\nabla^2 \Theta (X,X,Y,\theta Y)= - \nabla^2 \Theta (X,Y,X,\theta Y),$$
and, from Lemma~\ref{20} and Equation~\eqref{27},
$$ -\nabla^2 \Theta (X,Y,X,\theta Y) =-g(R(X,Y) X,Y)+ \eta(Y)g(\nabla^{2}_{X,X} \xi, Y)+ g(R(X,Y) \theta X, \theta Y),$$
therefore
\begin{align*}
&-g(R(X,Y) X,Y)+ \eta(Y)g(\nabla^{2}_{X,X} \xi, Y)+ g(R(X,Y) \theta X, \theta Y) \\
&=|(\nabla_X \theta)(Y)|^2 +\eta(Y) g(Y,\nabla^2_{X,X}\xi) +g^{2}(Y,\nabla_X \xi).
\end{align*}
\end{proof}

The somewhat cumbersome formula of the next proposition is key to the harmonicity of nearly cosymplectic structures, even though it is essentially nothing more than a refinement of Proposition~\ref{28}, followed by a fairly standard exercise on the symmetries of the Riemann curvature tensor.

\begin{proposition}\label{32}
Let $M^{2n+1}$ be a nearly cosymplectic manifold. Then for any $X,Y,Z$ and $W$ in $TM$, we have
$$R(W,X, Y, Z) - R(\theta W, \theta X, \theta Y, \theta Z) = \tfrac{1}{3}[A(W,X,Y,Z)-B(W,X,Y,Z)],$$
where
\begin{align*}
A(W,X,Y,Z)&=\tfrac{1}{2}[T(W+Y,Z+X)-T(W+Y,X)+T(W,X)+T(Y,X)\\
& -T(W,Z+X)-T(W+Y,Z)+T(W,Z)+T(Y,Z) - T(Y,Z+X)] ,
\end{align*}
\begin{align*}
B(W,X,Y,Z)&= \tfrac{1}{2} [ T(W+Z,X+Y) -T(W,X+Y ) - T(Z,X+Y)-T(W+Z,X) \\ 
& +T(W,X) +T(Z,X) -T(W+Z,Y) +T(W,Y)+T(Z,Y)]
\end{align*}
and
\begin{align*}
T(X,Y) &= -2g(Y,\xi)g(\theta (\nabla_X \theta) (Y), \nabla_X \xi)+ 2g(X,\xi)g(\theta (\nabla_X \theta) (Y), \nabla_Y \xi)\\
&-2g(X,\xi)g(Y,\xi)g(\nabla_X \xi , \nabla_Y \xi) + g^{2}(Y,\xi)| \nabla_X \xi|^2 + g^{2}(X,\xi)|\nabla_Y \xi|^2\\
& +g(Y,\xi)R(\theta X,\theta Y,X,\xi)-g(X,\xi) R(\theta X, \theta Y ,Y,\xi).
\end{align*}
\end{proposition}

\begin{proof}
From Lemma~\ref{nc7b}, we deduce that
\begin{align*} 
|(\nabla_{\theta X} \theta) (\theta Y)|^2 =& |(\nabla_{ X} \theta)(Y)|^2-2\eta(Y)g(\theta (\nabla_{X} \theta) (Y),\nabla_X \xi)\\
& +2 \eta(X) g( \theta (\nabla_{X} \theta)(Y),\nabla_Y \xi)\\
& -2\eta(X)\eta(Y)g(\nabla_X \xi , \nabla_Y \xi)\\
& + \eta^{2}(Y) |\nabla_X \xi|^2 + \eta^{2}(X) |\nabla_Y \xi|^2 ,
\end{align*}
since $\nabla_{\theta X} \xi = -\theta (\nabla_X \xi)$ and $\theta$ is skew-symmetric.
Substitute $|(\nabla_{X}\theta)(Y)|^2$ and $|(\nabla_{\theta X} \theta)(\theta Y)|^2$ using Proposition~\ref{28} to obtain
\begin{align*}
&-g(R(\theta X,\theta Y) \theta X, \theta Y ) +g(R(\theta X,\theta Y)  X,  Y)-g(R(\theta X, \theta Y) X, \eta (Y)\xi) \\
&-g(R(\theta X ,\theta Y) \eta(X)\xi, Y)-g^{2}(\nabla_X \xi, Y) \\ 
&=-g(R(X,Y) X, Y ) +g(R(X,Y) \theta X, \theta Y) -g^{2}(\nabla_X \xi, Y) -2\eta(Y)g(\theta (\nabla_{X} \theta) (Y),\nabla_X \xi)\\
& +2 \eta(X) g(\theta (\nabla_{X} \theta) (Y),\nabla_Y \xi)-2\eta(X)\eta(Y)g(\nabla_X \xi , \nabla_Y \xi)\\
& + \eta^{2}(Y) |\nabla_X \xi|^2 + \eta^{2}(X) |\nabla_Y \xi|^2 ,
\end{align*}
hence, for all $X$ and $Y$,
\begin{align*} 
&g(R(X,Y) X,Y)- g(R(\theta X, \theta Y) \theta X, \theta Y)= -2\eta(Y)g( \theta (\nabla_{X} \theta)(Y),\nabla_X \xi)\\
& +2 \eta(X) g((\theta \nabla_{X} \theta)(Y),\nabla_Y \xi)-2\eta(X)\eta(Y)g(\nabla_X \xi , \nabla_Y \xi)\\
& + \eta^{2}(Y) |\nabla_X \xi|^2 + \eta^{2}(X) | \nabla_Y \xi|^2 +\eta(Y)g(R(\theta X, \theta Y) X, \xi)-\eta(X)g(R(\theta X,\theta Y)Y,\xi) ,
\end{align*}
and a polarization argument for four-tensors yields the result.
\end{proof}

\begin{remark}\label{33} 
One can easily check that
\begin{itemize}
\item $\forall X,Y\in TM$, $T(X,Y)=T(Y,X)$.
\item $\forall X,Y\in \mathcal{F}$, $T(X,Y)=0$, therefore if $X,Y,Z,W\in \mathcal{F}$, 
$$R(\theta W,\theta X, \theta Y ,\theta Z) = R( W, X,  Y , Z).$$
\item $\forall Y\in \mathcal{F}$, $$T(\xi, Y) = -|\nabla_Y \xi|^2.$$
\item $\forall X\in \mathcal{F}$, 
$$T(\xi+X,Y)= - |\nabla_Y \xi|^2 +2g(\theta (\nabla_X \theta) (Y), \nabla_Y \xi)-g(R(\theta X, \theta Y) Y, \xi).$$
\end{itemize}
\end{remark}

Judicious selection of vectors for the previous proposition promptly gives the harmonic section equations.

\begin{proposition}
A nearly cosymplectic manifold satisfies the first harmonic section equation~\eqref{hse1}.
\end{proposition}

\begin{proof}
Recall that Proposition~\ref{ric*} showed that the first harmonic section equation vanishes if and only if 
$$\ricci^*(\theta X,\theta Y) = \ricci^*(X,Y), \quad \forall X,Y \in \mathcal{F}.$$
Then, by Remark~\ref{33}
\begin{align*}
\ricci^*(\theta X, \theta Y) &= g(R(\theta X , F_i)\theta F_i , \theta^2 Y)\\
&=g(R( \theta^2 X,\theta F_i) \theta^2 F_i, \theta^3 Y)\\
&=g(R( X,\theta F_i) \theta^2 F_i, \theta Y)\\
&=\ricci^*(X,Y),
\end{align*}
working with the orthonormal basis $\{\theta F_i\}_{i=1,\dots,2n}$.
\end{proof}

\begin{proposition} \label{38} 
The unit vector field $\xi$ of a nearly cosymplectic manifold $M^{2n+1}$ is always harmonic (as a unit vector field). 
Furthermore $$R^{\mathcal{F}} (F_i, \theta F_i)\xi =0,$$
where $\{F_i\}_{i=1,\dots,2n}$ is an orthonormal frame of $\mathcal F$, and the second harmonic section equation~\eqref{hse2} is satisfied.
\end{proposition}

\begin{proof}
We apply Proposition~\ref{32} to the vectors $W=F_i, X= \theta F_i, Y=\xi$ and $Z=\theta W$.
Then, using repeatedly Remark~\ref{33}, we compute the various terms to obtain
\begin{align*} 
&g(\theta (\nabla_{F_i} \theta) (\theta W) ,\nabla_{ \theta F_i}\xi ) + \tfrac{1}{2} g(R(\theta F_i, W) \theta F_i,\xi ) + \tfrac{1}{2} g(R(\theta E_i, E_i) \theta W,\xi)\\ 
&=g(\theta (\nabla_{\theta F_i} (\theta W)), \nabla_{F_i} \xi) +\tfrac{1}{2}g(R(F_i,\theta F_i) \theta W,\xi) -\tfrac{1}{2} g(R(F_i, W) F_i,\xi)\\
&+3g(R(F_i,\theta F_i)\xi ,\theta W).
\end{align*}
Since $g(R(F_i, W) F_i,\xi )=g(R(\theta F_i, W) \theta F_i,\xi )$, this simplifies to
\begin{align*} 
&g(\theta(\nabla_{F_i}\theta) (\theta W),\nabla_{\theta F_i}\xi )- g(\theta (\nabla_{\theta F_i} (\theta W)), \nabla_{F_i} \xi) 
- g(R(\xi, F_i) F_i,W ) \\ 
&= 2g(R(F_i,\theta F_i)\xi ,\theta W).
\end{align*}
Furthermore
\begin{align*} 
g(\theta (\nabla_{F_i} \theta) (\theta W) ,\nabla_{\theta F_i}\xi )&=-g((\nabla_{F_i} \theta) (\theta W) ,\nabla_{F_i}\xi )\\
&=-g(W,\theta (\tr (\nabla \theta \circ \nabla \xi))),
\end{align*}
and
\begin{align*} 
-g(\theta (\nabla_{\theta F_i} \theta) (\theta W) ,\nabla_{F_i}\xi )&=-g(-\theta^2 (\nabla_{\theta F_i} \theta)(W),\nabla_{F_i}\xi )\\
&=g((\nabla_{\theta F_i} \theta) (\theta^2W), \nabla_{F_i}\xi)\\
&=-g(W,\theta (\tr (\nabla \theta \circ \nabla \xi))),
\end{align*}
then the formula of Proposition~\ref{32} can be rewritten 
$$ 2g (W , \theta (\tr (\nabla \theta \circ \nabla \xi))) =2g(\theta R(F_i,\theta F_i)\xi, W )-g((\nabla^*\nabla \xi)^{\mathcal{F}},W).$$
On the other hand, in Proposition~\ref{nc10} we proved that
$$\theta \circ \tr ( \nabla \theta \circ  \nabla \xi) = -\tfrac{1}{2}\theta R(F_i, \theta F_i)\xi -(\nabla^*\nabla \xi)^{\mathcal{F}},$$
then
$$\theta R(F_i,\theta F_i)\xi-\tfrac{1}{2}(\nabla^*\nabla \xi)^{\mathcal{F}}=-\tfrac{1}{2} \theta R(F_i, \theta F_i)\xi -(\nabla^*\nabla \xi)^{\mathcal{F}},$$
that is
\begin{equation}\label{35}
3\theta R(F_i,\theta F_i)\xi = -(\nabla^*\nabla \xi)^{\mathcal{F}}.
\end{equation}
Now we use the formula of Proposition~\ref{32} with $W=\xi, X=F_i, Y=F_i$ and $Z\in \mathcal{F}$: 
$$g(R(\xi,F_i) F_i,Z)= -g(Z, \theta (\nabla_{F_i}\theta)(\nabla_{F_i} \xi))-\tfrac{1}{2}g(R(\theta F_i, \theta Z)F_i ,\xi).$$
Since we proved in Proposition~\ref{nc10} that
$$\theta\circ\tr(\nabla\theta\circ\nabla \xi) = -\tfrac{1}{2} \theta R(F_i, \theta F_i)\xi -(\nabla^*\nabla \xi)^{\mathcal{F}} ,$$
we obtain
$$ g(R(\xi,F_i) F_i,Z)=\tfrac{1}{2} g(\theta R(F_i, \theta F_i)\xi+2\nabla^*\nabla \xi,Z )+\tfrac{1}{2} g(R( \xi,F_i)\theta F_i, \theta Z),$$
and, since $\xi$ is a Killing vector field, 
\begin{equation}\label{37}
R^{\mathcal{F}} (F_i, \theta F_i)\xi = - R^{\mathcal{F}} (\xi,F_i) \theta F_i .
\end{equation}
By the first Bianchi identity and replacing the basis $\{F_i\}_{i=1,\dots,2n}$ by $\{\theta F_i\}_{i=1,\dots,2n}$ we have
$$ R(F_i, \theta F_i) \xi = -2R(\xi,F_i) \theta F_i .$$
However by Equation~\eqref{37} 
$$R^{\mathcal{F}} (F_i, \theta F_i)\xi = - R^{\mathcal{F}}(\xi,F_i) \theta F_i ,$$
and comparing these last two equations, we infer that
$$ R^{\mathcal{F}}(F_i,\theta F_i)\xi =0,$$
and, with Equation~\eqref{35}, conclude that
$$(\nabla^*\nabla \xi)^{\mathcal{F}}=0.$$
Moreover, from Proposition~\ref{nc10}, the second harmonic section equation must hold.
\end{proof}

\section{Harmonic map Equations}

Having just established that a nearly cosymplectic structure must be a harmonic section, the question naturally arises whether it is possible to go a little further and show that such an almost contact structure can actually be a harmonic map. Though the harmonic map equation, subordinated to the harmonic section equations, is given by the vanishing of a sum of two terms, we are able to prove that each summand vanishes, the first being entirely dependent on a complement to Proposition~\ref{32}, while the second hinges on the vertical part of the curvature tensor of $\mathcal F$.

\begin{proposition} \label{42}
Let $M^{2n+1}$ be a nearly cosymplectic manifold. Then for $X,Y,Z$ and $W$ in $\mathcal F$, we have
\begin{align*}
R(Y,X, W,Z) - R(Y,X, \theta W, \theta Z)&=-g((\nabla_W \theta) (Z) , (\nabla_Y \theta) (X)) \\
&+ g(Y,\nabla_X \xi) g(Z,\nabla_W \xi) .
\end{align*}
\end{proposition}

\begin{proof}
Let $X,Z$ and $W$ be vectors in the $\mathcal F$-distribution, applying Proposition~\ref{28} to the vectors $W$ and $Z+X$, we have
\begin{align*}
&|(\nabla_W \theta)(Z+X) |^2 +g^2(Z+X,\nabla_W \xi)\\
&=-g(R(W,Z+X)W,Z+X)+g(R(W,Z+X)\theta W,\theta (Z+X)). 
\end{align*}
Expanding both sides of the equation, we get
\begin{align*} 
& |(\nabla_W \theta) (Z)|^2 + |(\nabla_W \theta)(X)|^2 +2g((\nabla_W \theta)(Z),(\nabla_W \theta)(X) )\\
&+ g^2 (Z,\nabla_W \xi) +g^2 (X,\nabla_W \xi)+ 2 g(Z,\nabla_W \xi)g(X,\nabla_W \xi)\\
&= -g(R(W,Z)W,Z)+g(R(W,Z)\theta W , \theta Z)-2g(R(W,Z)W,X) \\
& +2g(R(W,Z)\theta W, \theta X) -g(R(W,X)W,X)+R(W,X,\theta W,\theta X) ,
\end{align*}
and since (Proposition~\ref{32})
$$g(R(W,X) \theta W , \theta Z) =g(R(W,Z) \theta W,\theta X),$$ 
we deduce that
\begin{align}
2g((\nabla_W \theta)(Z), (\nabla_W \theta) (X) ) + 2 g(Z,\nabla_W \xi)g(X,\nabla_W \xi)\label{eq2}\\
=-2g(R(W,Z) W,X) +2g(R(W,Z) \theta W, \theta X ).\notag
\end{align}

Applying this equation to the vectors $W, X+Y$ and $Z$, we obtain
\begin{align*} 
g(R(W, X+Y)(X+Y),Z)-g(R(W,X+Y) \theta (X+Y) ,\theta Z)\\
= g((\nabla_{X+Y} \theta) (W) ,(\nabla_{X+Y} \theta) (Z)) +g(W,\nabla_{X+Y} \xi)g(Z,\nabla_{X+Y} \xi), 
\end{align*}
and expanding the left-hand side, with the first Bianchi identity, we can rewrite this equation as
\begin{align*} 
&g(R(W,X) X,Z)-g(R(W,X)\theta  X, \theta Z)+ 2g(R(W,X)Y,Z)-g(R(W,X)\theta Y, \theta Z)\\
&+ g(R(W,Y)Y , Z)-g(R(W,Y)\theta X ,\theta Z)-g(R(Y,X) W,Z)-g(R(W,Y) \theta Y, \theta Z)\\
&=g((\nabla_X \theta) (W),(\nabla_X \theta) (Z)) + g(W,\nabla_X \xi)g(Z,\nabla_X \xi)\\
&+g((\nabla_Y \theta) (W), (\nabla_Y \theta) (Z) )+g(W,\nabla_Y \xi)g(Z,\nabla_Y \xi)\\
&+2 g(R(W,X)Y , Z)-g(R(W,Y)\theta X ,\theta Z)-g(R(Y,X) W,Z)-g(R(W,X) \theta Y, \theta Z).
\end{align*}
Expanding both sides of the equation and using Equation~\eqref{eq2}, yields, after some simplifications
\begin{align*} 
&2g(R(W,X) Y,Z)=
g((\nabla_X \theta) (W), (\nabla_Y \theta) (Z))+g((\nabla_Y \theta) (W) ,(\nabla_X \theta) (Z))\\
&+g(W,\nabla_X \xi)g(Z,\nabla_Y \xi ) +g (W,\nabla_Y \xi)g(Z,\nabla_X \xi)\\
&+g(R(W,X)\theta Y, \theta Z)+g(R(Y,X)W,Z)+g(R(W,Y)\theta X, \theta Z) .
\end{align*}
Now, by Proposition~\ref{28}
\begin{align}\label{diamondsuit}
&|(\nabla_{W+Y} \theta) (X+Z)|^2+g^{2}(X+Z,\nabla_{W+Y} \xi) \\
&=-g(R(W+Y,X+Z) (W+Y), X+Z )+g(R(W+Y,X+Z)\theta (W+Y) ,\theta (X+Z)),\notag
\end{align}
Expanding both sides of the equation and applying Proposition~\ref{28} and Equation~\eqref{eq2}, it simplifies to
\begin{align*} 
&2g((\nabla_W \theta)(X) , (\nabla_Y \theta) (X)) +2 g((\nabla_W \theta )(X) ,(\nabla_Y \theta) (Z))+2g((\nabla_W \theta) (Z) , (\nabla_Y \theta) (X)) \\
&+2 g((\nabla_W \theta) (Z) ,(\nabla_Y \theta) (Z))+2g(X,\nabla_W \xi)g(Z,\nabla_Y \xi)+ 2(X,\nabla_Y \xi)g(Z,\nabla_W \xi)\\
&=2g((\nabla_X \theta) (Y) , (\nabla_X \theta) (W)) +2 g((\nabla_Z \theta) (W) ,(\nabla_Z \theta) (Y))\\
&+g(R(W,X)\theta Y,\theta Z) -3g(R(Y,X)W, Z)+2g(R(W,Z)\theta Y, \theta X)\\
&-g(R(W,Y)\theta X, \theta Z)-g((\nabla_X \theta) (W), (\nabla_Y \theta) (Z))-g((\nabla_Y\theta) (W),(\nabla_X \theta) (Z))\\
&-g(W,\nabla_X \xi)g(Z,\nabla_Y \xi)-g(W,\nabla_Y \xi)g(Z,\nabla_X \xi) ,
\end{align*}
and further, if we use nearly cosympleticity and the Killing vector field properties of $\xi$:
\begin{align}\label{diamondsuitdiamondsuit}
&3g (R(Y,X) W,Z) -g(R(W,X)\theta Y, \theta Z)-2 g(R(W,Z) \theta Y, \theta X)+g(R(W,Y)\theta X, \theta Z)\\
&=-g((\nabla_Y \theta) (W),(\nabla_X \theta) (Z) )- 2g((\nabla_W \theta) (Z) , (\nabla_Y \theta) (X)) - g((\nabla_W \theta) (X) , (\nabla_Y \theta) (Z))\notag\\
&-2g(X, \nabla_Y \xi )g(Z,\nabla_W \xi)-g(X,\nabla_W \xi)g(Z,\nabla_Y \xi)-g(W,\nabla_Y \xi)g(Z,\nabla_X \xi) .\notag
\end{align}
Replace $W$ by $\theta W$ and $Z$ by $\theta Z$ in the above formula to get
\begin{align} \label{spadesuitspadesuit}
&3g (R(Y,X)\theta W, \theta Z ) +g(R(\theta W,X)\theta Y,  Z)-2 g(R(\theta W, \theta Z) \theta Y, \theta X)-g(R(\theta W,Y)\theta X,  Z)\\
&=-g((\nabla_Y \theta) (\theta W),(\nabla_X \theta) ( \theta Z) )- 2g((\nabla_{\theta W} \theta) (\theta Z) , (\nabla_Y \theta)(X)) \notag\\
&- g((\nabla_{\theta W} \theta) (X) , (\nabla_Y \theta)(\theta Z)) -2g(X, \nabla_Y \xi )g(\theta Z,\nabla_{\theta W} \xi)\notag\\
&-g(X,\nabla_{\theta W} \xi)g(\theta Z,\nabla_Y \xi)-g(\theta W,\nabla_Y \xi)g(\theta Z,\nabla_X \xi) .\notag
\end{align}
By Lemma~\ref{nca}
\begin{align*}
&-g((\nabla_Y \theta) (\theta W),(\nabla_X \theta) (\theta Z)) = -g(\theta (\nabla_Y \theta) (W) , \theta (\nabla_X \theta) (Z))-g(W,\nabla_Y \xi)g(Z,\nabla_X \xi) \\
&=g((\nabla_Y \theta) (W) , -(\nabla_X \theta)(Z) +g ((\nabla_X \theta) (Z), \xi)\xi)-g(W,\nabla_Y \xi)g(Z,\nabla_X \xi) \\
&=-g((\nabla_Y \theta) (W),(\nabla_X \theta)(Z))+ g ((\nabla_X \theta) (Z), \xi)g((\nabla_Y \theta) (W),\xi)-g(W,\nabla_Y \xi)g(Z,\nabla_X \xi)\\
&=-g((\nabla_Y \theta) (W),(\nabla_X \theta) (Z))+ g(Z,\theta (\nabla_X \xi))g(W, \theta (\nabla_ Y\xi))-g(W,\nabla_Y \xi)g(Z,\nabla_X \xi),
\end{align*}  
and similarly 
\begin{align*}
-g((\nabla_{\theta W} \theta) (X) , (\nabla_Y \theta) (\theta Z))&=g((\nabla_X \theta) (W),(\nabla_Y \theta)(Z))- 
g(W,\theta (\nabla_X \xi))g(Z, \theta (\nabla_ Y\xi))\\
& +g(W,\nabla_X \xi)g(Z,\nabla_Y \xi) .
\end{align*}  
Plugging in the last two equations in \eqref{spadesuitspadesuit} and then subtracting it to \eqref{diamondsuitdiamondsuit}, we obtain
\begin{align}\label{blacksquare}
&-4g((\nabla_W \theta) (Z) , (\nabla_Y \theta) (X)) -4 g(X, \nabla_Y \xi) g(Z,\nabla_W \xi)\\
&=5g(R(Y,X) W,Z) -g (R(Y, \theta X) \theta W,Z)-5g(R(W,Z) \theta Y,\theta X) +g(R(X,\theta Y)\theta W, Z).\notag
\end{align} 
Replacing $Z$ by $\theta Z$ and $Y$ by $\theta Y$, dividing by 5 and using Proposition~\ref{32}, gives
\begin{align}\label{blacksquareblacksquare}
&-\tfrac{4}{5}[g((\nabla_W \theta)(\theta Z) , (\nabla_{\theta Y} \theta)(X)) + g(X, \nabla_{\theta Y}\xi ) g(\theta Z,\nabla_W \xi)]\\
&=g(R( Y, \theta X) \theta W,Z) -\tfrac{1}{5}g (R(Y,  X)  W,Z)-g(R(X,\theta Y) \theta W,Z) -\tfrac{1}{5}g(R(X, Y)\theta W, \theta Z),\notag
\end{align} 
now adding up \eqref{blacksquare} and \eqref{blacksquareblacksquare} yields
\begin{align} \label{blacklozenge}
&-4g((\nabla_W \theta) (Z), (\nabla_Y \theta) (X))-\tfrac{4}{5}g((\nabla_W \theta) (\theta Z), (\nabla_{\theta Y} \theta) (X) )\\
&-4 g(X,\nabla_Y \xi)g(Z, \nabla_W \xi)- \tfrac{4}{5}g(X,\nabla_{\theta Y} \xi)g(\theta Z, \nabla_W \xi)\notag\\
&=\tfrac{24}{5}g(R(Y,X) W,Z) -\tfrac{24}{5}g(R(Y,X)\theta W , \theta Z).\notag
\end{align}
Finally, since 
\begin{align*}
g((\nabla_W \theta) (\theta Z), (\nabla_{\theta Y} \theta) (X) )&= -g((\nabla_W \theta) (Z), (\nabla_X \theta) (Y) )+ g(\theta Y , \nabla_X \xi)g(\theta Z , \nabla_W \xi)\\
&-g(Z,\nabla_W \xi) g(Y,\nabla_X \xi),
\end{align*}
injecting this term in Equation~\eqref{blacklozenge}, and simplifying by 5, yields the proposition.
\end{proof}

\begin{corollary} \label{43} 
For any $X,Y,Z$ and $W$ in $\mathcal{F}$ we have:
\begin{align*}
g([R(Y,X),\theta ] W, Z)&= \Theta ((\nabla_Y \theta) (X) , (\nabla_W \theta)(Z)) \\
&-g(W,\nabla_Z \xi)g(X, \nabla_{\theta Y}\xi)+g(Y,\nabla_X \xi)g(Z,\nabla_{\theta W} \xi) .
\end{align*}
\end{corollary}
\begin{proof} Substitute $W$ by $\theta W$ in Proposition~\ref{42} and use Lemma~\ref{nca} and the skew-symmetry of $\theta$ and $\nabla \theta$, to obtain
\begin{align*}
&g(R(Y,X) \theta W,Z) -g(\theta R(Y,X) W, Z)\\
&=-g(\theta (\nabla_{Z} \theta) (W),(\nabla_Y \theta) (X))+ g(W,\nabla_Z \xi)g( (\nabla_Y \theta) (X),\xi) + g(Y,\nabla_X \xi) g(Z,\nabla_{\theta W} \xi)\\
&= \Theta ((\nabla_Y \theta) (X),(\nabla_{W} \theta) (Z))- g(W,\nabla_Z \xi)g( X,\nabla_{\theta Y}\xi) + g(Y,\nabla_X \xi) g(Z,\nabla_{\theta W} \xi),
\end{align*}
which is the expression we sought.
\end{proof}

Recall from the introduction that the harmonic map equations for an almost contact structure are (\cite{VW1})
$$\langle (\bar{\nabla}_{E_i} J) (F_j), [\bar{R} (E_i, X), \theta ] F_j \rangle+ 8   \langle \nabla_{E_i} \xi, R(E_i,X)\xi \rangle =0 ,$$
for any $X\in TM$,and where $\{E_i\}_{i=1,\dots,2n+1}= \{F_i\}_{i=1,\dots,2n}\cup \{\xi\}$ an orthonormal basis of $TM$ with $\{F_i\}_{i=1,\dots,2n}$ an orthonormal basis of $\mathcal F$. 

We will actually prove that each term of this sum is zero, separately for $X$ in the ${\mathcal F}$-distribution and then for $X=\xi$ the Reeb vector field.

\begin{proposition} 
Let $M^{2n+1}$ be a nearly cosymplectic manifold. Then for any $X$ in $\mathcal{F}$ we have
\begin{align*} 
g( (\bar{\nabla}_{E_i} J) (F_j), [R(E_i, X), \theta ] F_j ) &= 0 ,
\end{align*}
where $\{F_i\}_{i=1,\dots,2n}$ is an orthonormal basis of $\mathcal F$ and $\{E_i\}_{i=1,\dots,2n+1}= \{F_i\}_{i=1,\dots,2n}\cup \{\xi\}$. 
\end{proposition}

\begin{proof}
First notice that if $X\in \mathcal{F}$ then, by Corollary~\ref{43}, we have
\begin{align*} 
&\langle (\bar{\nabla}_{E_i} J) (F_j), [R^{\mathcal{F}} (E_i, X), \theta ] F_j \rangle
=\Theta ((\nabla_{E_i} \theta) (X) , (\nabla_{F_j} \theta) ((\bar{\nabla}_{E_i} J) (F_j))) \\
&-g(F_j,(\nabla_{\bar{\nabla}_{E_i} J) (F_j)}\xi)g(X, \nabla_{\theta E_i}\xi)
+g(E_i,\nabla_X \xi)g(\bar{\nabla}_{E_i} J (F_j),\nabla_{\theta F_j} \xi)
\end{align*}
Now
\begin{align*} 
&g((\nabla_{F_i} \theta) (X) , \theta (\nabla_{F_j} \theta) \circ (\bar{\nabla}_{F_i}J)(F_j) )  
= g((\nabla_{X} \theta) (F_i) ,\theta (\nabla_{F_j} \theta) \circ (\bar{\nabla}_{F_j}J)(F_i))\\
& = -g(\theta (\nabla_{X} \theta) (F_i),(\nabla_{F_j} \theta )^2(F_i)) 
+ g((\nabla_{F_j} \theta) (F_i), \xi) g(\theta (\nabla_{X}\theta) (F_i),(\nabla_{F_j}\theta)(\xi)),
\end{align*}
and
$$ g(\theta (\nabla_{X} \theta) (F_i) , (\nabla_{F_j} \theta)^2(F_i) )=0,$$
since $(\nabla_{F_i} \theta )^2$ is a symmetric operator on $\mathcal{F}$, whilst $\theta \circ\nabla_ Y \theta$ is antisymmetric.
Furthermore, changing the basis $\{F_i\}_{i=1,\dots,2n}$ into the basis $\{\theta F_i\}_{i=1,\dots,2n}$ and using the Lemmas~\ref{nc7b} and \ref{nca} shows that
$$g((\nabla_{F_j} \theta) (F_i), \xi) g(\theta (\nabla_{X}\theta) (F_i),(\nabla_{F_j}\theta)(\xi)) =0.$$
On the other hand
\begin{align*} 
g((\nabla_{\xi} \theta) (X) , \theta (\nabla_{F_j} \theta) \circ (\bar{\nabla}_{\xi}J)(F_j)) & = 
-g(\nabla_{X} \xi ,(\nabla_{F_j} \theta) \circ (\bar{\nabla}_{\xi}J)(F_j))\\
&=-g( \nabla_{X} \xi ,(\nabla_{F_j} \theta) (\theta \nabla_{F_j} \xi))\\
&=-g( \nabla_{X} \xi , -\theta(\nabla_{F_j}\theta) (\nabla_{F_j} \xi) + g(\nabla_{F_j} \xi,\nabla_{F_j}\xi)\xi )\\
&=0,
\end{align*}
since $\xi$ is a unit section and (Proposition~\ref{38}) $\tr \nabla \theta \circ \nabla \xi =0$. 
Therefore, if $X\in \mathcal{F}$
$$g((\nabla_{E_i} \theta) (X) , \theta (\nabla_{F_j} \theta) \circ (\bar{\nabla}_{E_i}J)(F_j) )=0.$$
Our condition now reduces to
\begin{align*} 
\langle (\bar{\nabla}_{E_i} J) (F_j), [R^{\mathcal{F}} (E_i, X), \theta ] F_j \rangle
&= -g(F_j,\nabla_{(\bar{\nabla}_{E_i} J) (F_j)}\xi)g(X, \nabla_{\theta E_i}\xi)\\
&+g(E_i,\nabla_X \xi)g((\bar{\nabla}_{E_i} J)(F_j),\nabla_{\theta F_j} \xi),
\end{align*}
and concentrating on the second term
\begin{align*}
g(E_i,\nabla_{X} \xi )g((\bar{\nabla}_{E_i} J) (F_j),\nabla_{\theta F_j}\xi)&=
g(X, \nabla_{E_i} \xi)g(\theta F_j, \nabla_{(\bar{\nabla}_{E_i} J)(F_j)} \xi)\\
&= g(X, \nabla_{\theta E_i} \xi)g(\theta F_j, \nabla_{(\bar{\nabla}_{\theta E_i} J)(F_j)} \xi)\\
&+ g(X, \nabla_{\xi} \xi)g(\theta F_j, \nabla_{(\bar{\nabla}_{\xi} J)(F_j)} \xi) \\
& =g(X, \nabla_{\theta E_i} \xi)g(\theta F_j, \nabla_{\theta (\bar{\nabla}_{F_j} J)(E_i)} \xi)\\
& =-g(X, \nabla_{\theta E_i} \xi)g( F_j,  \nabla_{(\bar{\nabla}_{F_j} J)(E_i)} \xi)\\
& = g(X, \nabla_{\theta E_i} \xi)g( F_j,  \nabla_{(\bar{\nabla}_{E_i} J)(F_j)} \xi)
\end{align*}
and the result follows as it cancels with the first term. 
\end{proof}

We now prove the counterpart formula for vectors in the $\xi$-direction.

\begin{proposition}
Let $M^{2n+1}$ be a nearly cosymplectic manifold. Then 
$$ g([R(E_i,\xi) ,\theta ] F_j, (\bar{\nabla}_{E_i} J)(F_j)) =0,$$
where $\{F_i\}_{i=1,\dots,2n}$ is an orthonormal basis of $\mathcal F$ and $\{E_i\}_{i=1,\dots,2n+1}= \{F_i\}_{i=1,\dots,2n}\cup \{\xi\}$. 
\end{proposition}

\begin{proof}If $W$ and $Z$ are in the $\mathcal{F}$-distribution, we have 
$$ g([R(W,\xi) ,\theta ] Y, Z) =g(R(W,\xi) \theta Y,Z) + g( R (W,\xi) Y,\theta Z ),$$
and using Proposition~\ref{32} with $W$, $Y$ and $Z$ in $\mathcal F$ and $\xi$, we obtain
\begin{align} 
&3g(R(W,\xi) Y,Z) = 2g(\theta (\nabla_{Z} \theta) (Y) ,\nabla_{W} \xi)+ g(\theta (\nabla_{Z} \theta) (W) ,\nabla_{Y} \xi) \label{eqp45}\\
& -\tfrac{1}{2}g(R(\theta Z, \theta W) Y, \xi) + g (R(\theta Y, \theta Z)W , \xi)
-g(\theta (\nabla_{Y} \theta) (W) ,\nabla_{Z} \xi)+\tfrac{1}{2}g (R(\theta Y, \theta W)Z , \xi) \notag .
\end{align}
Putting $W=E_i$, $Y= \theta F_j$ and $Z=(\bar{\nabla}_{E_i} J) (F_j)$ and then $W=E_i$, $Y= F_j$ and $Z=\theta (\bar{\nabla}_{E_i} J)(F_j)$, we obtain
\begin{align*} 
&3g(R(E_i ,\xi) \theta F_j,(\bar{\nabla}_{E_i} J) (F_j)) = \\
&2g(\theta (\nabla_{(\bar{\nabla}_{E_i} J)(F_j)} \theta) (\theta F_j) ,\nabla_{E_i} \xi)+ 
g(\theta (\nabla_{(\bar{\nabla}_{E_i} J) (F_j)} \theta) (E_i) ,\nabla_{\theta F_j} \xi)\\
& -\tfrac{1}{2}g(R(\theta (\bar{\nabla}_{E_i} J) (F_j), \theta E_i) \theta F_j, \xi) - g (R( F_j, \theta (\bar{\nabla}_{E_i} J) (F_j))E_i , \xi)\\
&-g(\theta (\nabla_{\theta F_j} \theta) (E_i) ,\nabla_{(\bar{\nabla}_{E_i} J) (F_j)} \xi)-\tfrac{1}{2}g(R(F_j, \theta E_i)(\bar{\nabla}_{E_i} J) (F_j), \xi) 
\end{align*}
and 
\begin{align*} 
&3g(R(E_i,\xi)  F_j, \theta (\bar{\nabla}_{E_i} J) (F_j)) = \\
&2g(\theta (\nabla_{\theta (\bar{\nabla}_{E_i} J) (F_j)} \theta) ( F_j) ,\nabla_{E_i} \xi)+ g(\theta (\nabla_{\theta (\bar{\nabla}_{E_i} J) (F_j)} \theta) (E_i) ,\nabla_{ F_j} \xi)\\
& +\tfrac{1}{2}g(R( (\bar{\nabla}_{E_i} J) (F_j), \theta E_i)  F_j, \xi) - g (R(\theta F_j,  (\bar{\nabla}_{E_i} J) (F_j))E_i , \xi)\\
&-g(\theta (\nabla_{ F_j} \theta) (E_i) ,\nabla_{\theta (\bar{\nabla}_{E_i} J) (F_j)} \xi)+\tfrac{1}{2}g (R(\theta F_j, \theta E_i)\theta (\bar{\nabla}_{E_i} J) (F_j), \xi) .
\end{align*}
The first terms of each equation are opposite, since
\begin{align*}
2g(\theta (\nabla_{ (\bar{\nabla}_{E_i} J) (\theta F_j)} \theta) ( F_j) ,\nabla_{E_i} \xi) &= 
2g( (\nabla_{ (\bar{\nabla}_{E_i} J) (F_j)} \theta) ( F_j) ,\nabla_{E_i} \xi)
\end{align*}
and
\begin{align*}
2g(\theta (\nabla_{\theta (\bar{\nabla}_{E_i} J) (F_j)} \theta) ( F_j),\nabla_{E_i} \xi)&=
-2g(\theta^2 (\nabla_{ (\bar{\nabla}_{E_i} J) ( F_j)} \theta) ( F_j) ,\nabla_{E_i} \xi)\\
&= 2g( (\nabla_{(\bar{\nabla}_{E_i} J) ( F_j)} \theta_ ( F_j) , \nabla_{E_i} \xi),
\end{align*}
by Lemmas~\ref{nca} and \ref{nc7b}.

The second terms will cancel since
\begin{align*}
g(\theta (\nabla_{(\bar{\nabla}_{E_i} J) (F_j)} \theta) (E_i), \nabla_{\theta F_j} \xi) &=
- g((\nabla_{(\bar{\nabla}_{E_i} J) (F_j)} \theta) (E_i),  \nabla_{ F_j} \xi) ,
\end{align*}
and
\begin{align*} 
g(\theta (\nabla_{\theta (\bar{\nabla}_{E_i} J) (F_j)} \theta) (E_i), \nabla_{ F_j} \xi) &= 
g( (\nabla_{ (\bar{\nabla}_{E_i} J) (F_j)} \theta) ( E_i), \nabla_{ F_j} \xi),
\end{align*}
as will the fifth terms
\begin{align*} 
-g(\theta (\nabla_{ \theta F_j} \theta) (E_i) ,\nabla_{(\bar{\nabla}_{E_i} J) (F_j)} \xi) &=  
-g( (\nabla_{F_j} \theta) (E_i) ,\nabla_{(\bar{\nabla}_{E_i} J) (F_j)} \xi)
\end{align*}
whilst
\begin{align*} 
-g(\theta (\nabla_{F_j} \theta) (E_i) ,\nabla_{\theta (\bar{\nabla}_{E_i} J) (F_j)} \xi) &= 
g( (\nabla_{F_j} \theta) (E_i) , \nabla_{(\bar{\nabla}_{E_i} J) (F_j)} \xi).
\end{align*}
In conclusion we are only left with the curvature terms
\begin{align*}
&3g([R(E_i, \xi) , \theta ] F_j , (\bar{\nabla}_{E_i} J)(F_j))\\
&=-\tfrac{1}{2}g(R(\theta (\bar{\nabla}_{E_i} J)(F_j), \theta E_i)\theta F_j,\xi)-\tfrac{1}{2}g(R(F_j, \theta E_i) (\bar{\nabla}_{E_i} J)(F_j),\xi)\\
&+\tfrac{1}{2} g(R(\theta F_j,\theta E_i) \theta (\bar{\nabla}_{E_i} J)(F_j), \xi)+\tfrac{1}{2}g(R((\bar{\nabla}_{E_i} J)(F_j), \theta  E_i) F_j,\xi )\\
&-g(R(F_j, \theta (\bar{\nabla}_{E_i} J)(F_j) ) E_i, \xi)-g(R(\theta F_j , (\bar{\nabla}_{E_i} J)(F_j))E_i ,\xi).
\end{align*}
Using the first Bianchi identity this simplifies to
\begin{align}\label{surd}
&3g([R(E_i, \xi) , \theta ] F_j , (\bar{\nabla}_{E_i} J)(F_j))\\
&=\tfrac{1}{2}g(R((\bar{\nabla}_{E_i} J)(F_j),F_j) \theta E_i,\xi)
+\tfrac{1}{2} g(R(\theta F_j, \theta (\bar{\nabla}_{E_i} J)(F_j) ) \theta E_i , \xi)\notag\\
&-g(R(F_j, \theta (\bar{\nabla}_{E_i} J)(F_j) ) E_i, \xi)-g(R(\theta F_j , (\bar{\nabla}_{E_i} J)(F_j))E_i ,\xi)\notag .
\end{align}
Changing the  orthonormal basis $\{E_i\}_{i=1,\dots,2n+1}$ to $\{\theta E_i \}_{i=1,\dots,2n}\cup \{\xi\}$ yields
\begin{align*} 
\tfrac{1}{2}g(R((\bar{\nabla}_{E_i} J)(F_j),F_j) \theta E_i,\xi)
&=\tfrac{1}{2}g(R((\bar{\nabla}_{\theta E_i} J)(\theta F_j),\theta F_j) \theta^2 E_i,\xi)\\
&=-\tfrac{1}{2}g(R(\theta F_j, (\bar{\nabla}_{E_i} J)(F_j)) E_i,\xi) ,
\end{align*}
and making the same type of change
\begin{align*}
\tfrac{1}{2} g(R(\theta F_j, \theta (\bar{\nabla}_{E_i} J)(F_j) ) \theta E_i , \xi) &= 
-\tfrac{1}{2} g(R( F_j, \theta (\bar{\nabla}_{E_i} J)( F_j) )  E_i , \xi),
\end{align*}
therefore Equation~\eqref{surd} reduces to
\begin{align}\label{surdsurd}
&3g([R(E_i, \xi) , \theta ] F_j , (\bar{\nabla}_{E_i} J)(F_j))\\
&=-\tfrac{3}{2} g(R(F_j, \theta (\bar{\nabla}_{E_i} J)(F_j) ) E_i, \xi)-\tfrac{3}{2} g(R(\theta F_j , (\bar{\nabla}_{E_i} J)(F_j))E_i ,\xi).\notag
\end{align}
Changing from $\{F_i\}_{i=1,\dots,2n}$ to $\{\theta F_i\}_{i=1,\dots,2n}$, we can rewrite the second term:
\begin{align*} 
-\tfrac{3}{2} g(R(\theta F_j , (\bar{\nabla}_{E_i} J)(F_j))E_i ,\xi)&= -\tfrac{3}{2} g(R( F_j ,\theta(\bar{\nabla}_{E_i} J)(F_j))E_i ,\xi),
\end{align*}
and Equation~\eqref{surdsurd} becomes
\begin{align}\label{surdsurdsurd}
g([R(E_i, \xi) , \theta ] F_j , (\bar{\nabla}_{E_i} J)(F_j))= -g(R( F_j ,\theta (\bar{\nabla}_{E_i} J)(F_j))E_i ,\xi).
\end{align}
Now
\begin{align}\label{wp}
&g([R(E_i, \xi) , \theta ] F_j , (\bar{\nabla}_{E_i} J)(F_j))= 
g(R(E_i, \xi) \theta F_j, (\bar{\nabla}_{E_i} J) (F_j)) - g(\theta R(E_i, \xi)  F_j, (\bar{\nabla}_{E_i} J) (F_j))\\
&=g(R(E_i, \xi) \theta^2 F_j, (\bar{\nabla}_{E_i}J)(\theta F_j))+g( R(E_i, \xi)F_j, \theta (\bar{\nabla}_{E_i} J) ( F_j))\notag\\
&= 2g(R(E_i, \xi) F_j, \theta (\bar{\nabla}_{E_i} J) ( F_j)), \notag
\end{align}
again replacing $\{F_i\}_{i=1,\dots,2n}$ by $\{\theta F_i\}_{i=1,\dots,2n}$. Finally, comparing \eqref{wp} and \eqref{surdsurdsurd} implies that
$$g([R(E_i, \xi) , \theta ] F_j , (\bar{\nabla}_{E_i} J)(F_j)) = 0.$$
\end{proof}

We now prove that the second term in the harmonic map equation~\eqref{hme} vanishes. This requires a formula similar to Equation~\eqref{eqp45} but without curvature terms on the right-hand side, and the proposition will easily follow from this.

\begin{proposition} 
Let $M^{2n+1}$ be a nearly cosymplectic manifold. Then for any $X$ in $TM$, we have
$$g ( \nabla_{E_i} \xi , R(E_i, X) \xi) =0 ,$$
where $\{E_i\}_{i=1,\dots,2n+1}= \{F_i\}_{i=1,\dots,2n}\cup \{\xi\}$ is an orthonormal basis of $TM$ with $\{F_i\}_{i=1,\dots,2n}$ an orthonormal basis of $\mathcal F$.
\end{proposition}

\begin{proof}
Let $Y,Z$ and $W$ be vectors in the $\mathcal F$-distribution.
We first need an expression for $g(R(W,Z)Y,\xi)$.
Apply Proposition~\ref{28} to $W+Y$ and $Z + \xi$:
\begin{align*} 
&-g(R(W+Y, \xi+Z)(W+Y),\xi+Z)+g(R(W+Y, \xi+Z)\theta (W+Y),\theta Z)\\
&= |(\nabla_{W+Y} \theta)(\xi+Z)|^2+g^2(\xi+Z ,\nabla_{W+Y} \xi) ,
\end{align*}
Expanding the right-hand side of the equation, with the help of the first Bianchi identity, we obtain
\begin{align}\label{@}
&-g(R(W+Y, \xi+Z)(W+Y),\xi+Z)+g(R(W+Y, \xi+Z)\theta (W+Y),\theta Z)\\
&= -g(R(W,\xi) W,\xi)-g(R(W,Z) W,Z )+g(R(W,Z) \theta W ,\theta Z) -g(R(Y,\xi) Y,\xi) \notag\\
&-g(R(Y,Z) Y,Z )+g(R(Y,Z) \theta Y ,\theta Z)-2g(R(W,\xi) W,Z)+g(R(W,\xi) \theta W,\theta Z )\notag\\
&-2g(R(W,Z) Y,Z)+2g(R(W,Z) \theta Y,\theta Z)-2g(R(Y,\xi) Y,Z)+g(R(Y,\xi) \theta Y,\theta Z )\notag\\
&-2g(R(W,\xi) Y,\xi)-2g(R(W,\xi) Y,Z)+g(R(W,\xi) \theta Y,\theta Z )-2g(R(W,Z) Y,\xi)\notag\\
&+g(R(Y,\xi) \theta W,\theta Z ).\notag
\end{align}
But applying Proposition~\ref{28} to the vectors $W$ and $Z+\xi$ and expanding both sides we deduce that
\begin{equation}\label{frakA}
-2g((\nabla_W \theta) (Z) ,\theta \nabla_W  \xi) = -2g(R(W,\xi) W,Z) +g(R(W,\xi) \theta W, \theta Z),
\end{equation}
and similarly, with the vectors $W+Y$ and $\xi$, we obtain
\begin{equation}\label{frakB}
-g( R(\xi, W) \xi , Y) =g(\nabla_W \xi, \nabla_Y \xi).
\end{equation}
Therefore Equation~\eqref{@} becomes
\begin{align}\label{@@}
&-g(R(W+Y, \xi+Z)(W+Y),\xi+Z)+g(R(W+Y, \xi+Z)\theta (W+Y),\theta Z)\\
&= |(\nabla_W \theta)(\xi)|^2 + |(\nabla_Y \theta) (\xi)|^2+ |(\nabla_W \theta) (Z)|^2 +g^2(Z,\nabla_W \xi)
+|(\nabla_Y \theta) (Z)|^2 + g^2(Z, \nabla_Y \xi)\notag\\
&-2g((\nabla_W \theta) (Z) , \theta (\nabla_W \xi))+2 g(\nabla_W \xi, \nabla_Y \xi)-2g((\nabla_Y \theta) (Z) , \theta (\nabla_Y \xi))\notag\\
&+2 g((\nabla_Z \theta) (Y) , \nabla_Z \theta (W)) + 2g(Y,\nabla_Z \xi)g(W,\nabla_Z \xi) \notag\\
&-2g(R(W, \xi) Y ,Z )+g (R(W,\xi)\theta Y,\theta Z) -2g(R(Y,\xi) W,Z)+g(R(Y,\xi) \theta W, \theta Z).\notag
\end{align}
Since
\begin{align} \label{frakE}
&g(R(W,\xi+ Y) (\xi+Y) ,Z) -g(R(W,\xi+Y) \theta Y, \theta Z) \\
&=g(R(W,\xi) \xi, Z)+ 2g(R(W,\xi)Y,Z) -g(R(W,\xi)\theta Y, \theta Z)\notag\\
&-g(R(Y,W)Y, Z)+ g(R(Y,W)\theta Y, \theta Z)-g(R(Y,\xi) W,Z), \notag 
\end{align}
and by Proposition~\ref{42}
$$-g(R(Y,W)Y, Z)+ g(R(Y,W)\theta Y, \theta Z) = g((\nabla_Y \theta) (W) , (\nabla_Y \theta) (Z)) +g(W, \nabla_Y \xi) g(Z,\nabla_Y \xi), $$
we get
\begin{align} 
&g(R(W,\xi+ Y) (\xi+Y) ,Z) -g(R(W,\xi+Y) \theta Y, \theta Z) \label{skull} \\
&=g(\nabla_W \xi, \nabla_Z \xi)+2g(R(W,\xi) Y, Z) -g(R(W,\xi) \theta Y, \theta Z)\notag\\
&-g(R(Y,\xi) W, Z)+ g((\nabla_Y \theta) (W), (\nabla_Y \theta) (Z) )+ g(W,\nabla_Y \xi)g(Z,\nabla_Y \xi) . \notag
\end{align}
In order to compute the term~\eqref{skull}, we consider the vectors $W+Z$ and $Y+\xi$:
\begin{align}\label{maltese}
&g(R(W+Z) ,\xi+Y)(\xi +Y ),W+Z)- g(R(W+Z, \xi + Y) \theta Y, \theta (W+Z))\\
&=g(R(W,\xi +Y) (\xi+Y) ,W ) -g(R(W,\xi+Y) \theta Y , \theta W)\notag\\
&+2 g(R(W,\xi+Y ) (\xi+Y) , Z) -g (R(W , \xi +Y) \theta Y, \theta Z)\notag\\
&-g(R(Z,\xi +Y) \theta Y, \theta W )+g(R(Z,\xi+Y) (\xi+Y) ,Z) -g(R(Z, \xi+Y) \theta Y, \theta Z),\notag
\end{align}
and
\begin{align}
&-3g(R(W,\xi + Y ) \theta Y , \theta Z) +3 g(R(\theta W, \theta Y) Y,Z)\label{frakg} \\
&+3g(R(Z,\xi + Y ) \theta Y , \theta W) -3 g(R(\theta Z, \theta Y) Y,W) \label{frakh} \\
&= -3g(R(W,\xi+Y) \theta Y, \theta Z) +3 g(R(Z,\xi + Y)\theta Y , \theta W) .\notag 
\end{align}
Compute \eqref{frakg}, by Proposition~\ref{32} with $W$, $\xi + Y$, $\theta Y$ and $\theta Z$
\begin{align*} 
&3g(R(W,\xi+Y ) \theta Y, \theta Z) - 3 g(R(\theta W,\theta Y) Y,Z) \\
&= \tfrac{1}{2}[2g(\theta (\nabla_{\theta Z} \theta) (W), \nabla_{\theta Y} \xi)+ 4g(\theta (\nabla_{\theta Z} \theta) (\theta Y) , \nabla_{W} \xi)
-2g(\theta (\nabla_{\theta Y} \theta) (W) , \nabla_{\theta Z} \xi)\\
&+g(R(Z,\theta W) \theta Y, \xi)-2g(R(Z,Y) W,\xi) - g(R(Y,\theta W) \theta Z ,\xi)] ,
\end{align*}
and, exchanging $W$ and $Z$, we get \eqref{frakh}
\begin{align*} 
&3g(R(Z,\xi+Y ) \theta Y, \theta W) - 3 g(R(\theta Z,\theta Y) Y,W) \\
&= \tfrac{1}{2}[2g(\theta (\nabla_{\theta W} \theta) (Z), \nabla_{\theta Y} \xi)+ 4g(\theta (\nabla_{\theta W} \theta) (\theta Y) , \nabla_{Z} \xi)
-2g(\theta (\nabla_{\theta Y} \theta) (Z) , \nabla_{\theta W} \xi)\\
&+g(R(W,\theta Z) \theta Y, \xi)-2g(R(W,Y) Z,\xi) - g(R(Y,\theta Z) \theta W ,\xi)] ,
\end{align*}
hence
\begin{align}
&3g(R(W,\xi + Y ) \theta Y , \theta Z) - 3 g(R(\theta W, \theta Y) Y,Z) \label{davidsstar}\\
&-3g(R(Z,\xi + Y ) \theta Y , \theta W) + 3 g(R(\theta Z, \theta Y) Y,W) \notag\\
&= 3g(R(W,\xi+Y) \theta Y, \theta Z) -3 g(R(Z,\xi + Y)\theta Y , \theta W)\notag\\
&=\tfrac{1}{2}[2g(\theta (\nabla_{\theta Z} \theta) (W) , \nabla_{\theta Y} \xi)+4g(\theta (\nabla_{\theta Z} \theta) (\theta Y) , \nabla_W \xi) \notag\\
&+g(R(Z,\theta W) \theta Y, \xi) -2 g(\theta (\nabla_{\theta Y} \theta) (W) ,\nabla_{\theta Z} \xi) -g(R(Y, \theta W) \theta Z , \xi) -2g( R(Z,Y) W,\xi) \notag\\
&-2g(\theta (\nabla_{\theta W} \theta) (Z) , \nabla_{\theta Y} \xi)-4g(\theta (\nabla_{\theta W} \theta) (\theta Y), \nabla_Z \xi)\notag\\
&-g(R(W,\theta Z) \theta Y,\xi)+2g(\theta (\nabla_{\theta Y} \theta) (Z) , \nabla_{\theta W} \xi)+g(R(Y, \theta Z) \theta W, \xi) +2g(R(W, Y) Z, \xi)].\notag 
\end{align}
But
\begin{align*}
g(\theta (\nabla_{\theta Z} \theta) (W) , \nabla_{\theta Y} \xi) 
&= g((\nabla_{\theta W} \theta) (\theta Z) , \nabla_{\theta Y} \xi)\\
& = -g(\theta (\nabla_{\theta W} \theta) (Z), \nabla_{\theta Y} \xi),
\end{align*}
and
\begin{align*}
g(\theta (\nabla_{\theta Y} \theta) (Z) , \nabla_{\theta W} \xi)
&=-g((\nabla_{\theta Y} \theta)(Z),\nabla_W \xi)\\
&=-g((\nabla_{ Y} \theta) (\theta Z),\nabla_W \xi)\\
&=g((\nabla_{ \theta Z} \theta)(Y) ,  \nabla_W \xi)\\
&= g(\theta (\nabla_{\theta Z} \theta) (\theta Y) ,\nabla_W \xi).
\end{align*}
Therefore, we rewrite \eqref{davidsstar} as
\begin{align}\label{mercury}
&3g(R(W,\xi+Y) \theta Y, \theta Z) -3 g(R(Z,\xi + Y)\theta Y , \theta W)\\
&=\tfrac{1}{2}[4g(\theta (\nabla_{\theta Z} \theta) (W), \nabla_{\theta Y} \xi)+ 6g(\theta (\nabla_{\theta Z} \theta) (\theta Y) , \nabla_{W} \xi)
- 6g(\theta (\nabla_{\theta W} \theta) (\theta Y) , \nabla_{Z} \xi) \notag\\
&+g(R(Z,\theta W) \theta Y , \xi) -g(R(Y, \theta W) \theta Z , \xi)+2g(R(W,Y) Z,\xi)\notag\\
&-g(R(W,\theta Z) \theta Y, \xi) +g (R(Y,\theta Z) \theta W ,\xi) -2g(R(Z,Y) W,\xi)],\notag
\end{align}
and, using this equation, \eqref{maltese} becomes
\begin{align*}
&g(R(W+Z, \xi +Y ) (\xi+Y) ,W+Z) -g(R(W+Z, \xi +Y ) \theta Y ,\theta (W+Z)) \\
&= g(R(W,\xi+Y) (\xi+Y) ,W) -  g(R(W,\xi+Y) \theta Y ,\theta W) \\
&+2g(R(W,\xi+Y) (\xi+Y) ,Z) -2g(R(W,\xi+Y) \theta Y ,\theta Z)\\
&+\tfrac{1}{6}[4g(\theta (\nabla_{\theta Z} \theta) (W), \nabla_{\theta Y} \xi)+ 6g(\theta (\nabla_{\theta Z} \theta) (\theta Y) , \nabla_{W} \xi)
- 6g(\theta (\nabla_{\theta W} \theta) (\theta Y) , \nabla_{Z} \xi) \\
&+g(R(Z,\theta W) \theta Y , \xi) -g(R(Y, \theta W) \theta Z , \xi)+2g(R(W,Y) Z,\xi)\\
&-g(R(W,\theta Z) \theta Y, \xi) +g (R(Y,\theta Z) \theta W ,\xi) -2g(R(Z,Y) W,\xi)]\\
&+g(R(Z,\xi+Y) (\xi +Y) ,Z) -g(R(Z,\xi+Y) \theta Y , \theta Z)
\end{align*}
Using Proposition~\ref{28} and simplifying we obtain
\begin{align*}
&g(R(W,\xi+Y) (\xi+Y) ,Z) -g(R(W,\xi+Y) \theta Y ,\theta Z)=\\
&g(\nabla_{W} \xi  , \nabla_{Z} \xi)-g(\theta \nabla_W \xi ,(\nabla_Z \theta) (Y))-g((\nabla_W \theta) (Y) , \theta \nabla_Z \xi) \\
&+g((\nabla_W \theta) (Y) ,(\nabla_Z\theta)(Y))+g(Y,\nabla_W \xi) g(\nabla_Z \xi , Y)\\
&-\tfrac{1}{12}[4g(\theta (\nabla_{\theta Z} \theta) (W), \nabla_{\theta Y} \xi)+ 6g(\theta (\nabla_{\theta Z}\theta)(\theta Y),\nabla_{W} \xi)
- 6g(\theta (\nabla_{\theta W} \theta) (\theta Y) , \nabla_{Z} \xi) \\
&+g(R(Z,\theta W) \theta Y , \xi) -g(R(Y, \theta W) \theta Z , \xi)+2g(R(W,Y) Z,\xi)\\
&-g(R(W,\theta Z) \theta Y, \xi) +g (R(Y,\theta Z) \theta W ,\xi) -2g(R(Z,Y) W,\xi)].
\end{align*}
Plugging this into \eqref{skull} yields, after simplifications
\begin{align*} 
&2g(R(W,\xi) Y, Z) -g(R(W,\xi) \theta Y, \theta Z) \\
&= g(R(Y,\xi) W, Z)-g(\theta \nabla_W \xi ,(\nabla_Z \theta) (Y)) -g((\nabla_W \theta) (Y) , \theta \nabla_Z \xi) \\
&-\tfrac{1}{12}[4g(\theta (\nabla_{\theta Z} \theta) (W), \nabla_{\theta Y} \xi)+ 6g(\theta (\nabla_{\theta Z} \theta) (\theta Y) , \nabla_{W} \xi)\\
&- 6g(\theta (\nabla_{\theta W} \theta) (\theta Y) , \nabla_{Z} \xi) +g(R(Z,\theta W) \theta Y , \xi) -g(R(Y, \theta W) \theta Z , \xi)\\
&+2g(R(W,Y) Z,\xi)-g(R(W,\theta Z) \theta Y, \xi) +g (R(Y,\theta Z) \theta W ,\xi) -2g(R(Z,Y) W,\xi)].
\end{align*}
Then Equation~\eqref{@@} becomes
\begin{align*}
&-g(R(W+Y, \xi+Z)(W+Y),\xi+Z)+g(R(W+Y, \xi+Z)\theta (W+Y),\theta Z)\\
&= |(\nabla_W \theta) (\xi)|^2 + |(\nabla_Y \theta) (\xi)|^2 + |(\nabla_W \theta) (Z)|^2 +g^{2}(Z,\nabla_W \xi)+|(\nabla_Y \theta) (Z)|^2 \\
&+ g^{2}(Z, \nabla_Y \xi) -2g((\nabla_W \theta) (Z) , \theta (\nabla_W \xi))
+2 g(\nabla_W \xi, \nabla_Y \xi)-2g((\nabla_Y \theta)(Z) , \theta(\nabla_Y \xi))\\
&+2 g((\nabla_Z \theta)(Y) , (\nabla_Z \theta)(W)) + 2g(Y,\nabla_Z \xi)g(W,\nabla_Z \xi) \\
&-g(R(Y,\xi) W, Z)+g(\theta \nabla_W \xi ,(\nabla_Z \theta)(Y)) +g((\nabla_W \theta)(Y) , \theta \nabla_Z \xi) \\
&+\tfrac{1}{12}[4g(\theta (\nabla_{\theta Z}\theta)(W), \nabla_{\theta Y} \xi)+ 6g(\theta (\nabla_{\theta Z} \theta)(\theta Y) , \nabla_{W} \xi)\\
&- 6g(\theta (\nabla_{\theta W} \theta)(\theta Y) , \nabla_{Z} \xi) +g(R(Z,\theta W) \theta Y , \xi) -g(R(Y, \theta W) \theta Z , \xi)\\
&+2g(R(W,Y) Z,\xi)-g(R(W,\theta Z) \theta Y, \xi) +g (R(Y,\theta Z) \theta W ,\xi) -2g(R(Z,Y) W,\xi)]\\
&-2g(R(Y,\xi) W,Z)+g(R(Y,\xi) \theta W, \theta Z),
\end{align*}
and since
\begin{align*}
&-g(R(W+Y, \xi+Z)(W+Y),\xi+Z)+g(R(W+Y, \xi+Z)\theta (W+Y),\theta Z) \\
&= |(\nabla_{W+Y} \theta)(Z+ \xi)|^2 + g^{2}(W+Y,\nabla_{Z+\xi} \xi) ,
\end{align*}
we can simplify some terms to get
\begin{align*}
&2g((\nabla_W \theta)(Z) , \nabla_Y \theta (\xi)) + 2g((\nabla_Y \theta)(Z) , \nabla_W \theta (\xi))\\
&= -g(R(Y,\xi) W, Z)+g(\theta \nabla_W \xi ,(\nabla_Z \theta)(Y)) +g((\nabla_W \theta)(Y) , \theta \nabla_Z \xi) \\
&+\tfrac{1}{12}[4g(\theta (\nabla_{\theta Z}\theta)(W), \nabla_{\theta Y} \xi)+ 6g(\theta(\nabla_{\theta Z}\theta)(\theta Y),\nabla_{W} \xi)\\
&- 6g(\theta(\nabla_{\theta W}\theta)(\theta Y) , \nabla_{Z} \xi) +g(R(Z,\theta W) \theta Y , \xi) -g(R(Y, \theta W) \theta Z , \xi)\\
&+2g(R(W,Y) Z,\xi)-g(R(W,\theta Z) \theta Y, \xi) +g (R(Y,\theta Z) \theta W ,\xi) -2g(R(Z,Y) W,\xi)]\\
&-2g(R(Y,\xi) W,Z)+g(R(Y,\xi) \theta W, \theta Z).
\end{align*}
Finally, Equation~\eqref{@@} becomes
\begin{align*}
&2g((\nabla_W \theta)(Z),(\nabla_Y \theta)(\xi)) + 2g((\nabla_Y \theta)(Z) , (\nabla_W \theta)(\xi))\\
&=g(\theta \nabla_W \xi ,(\nabla_Z \theta)(Y)) +g((\nabla_W \theta)(Y) , \theta \nabla_Z \xi) \\
&+\tfrac{1}{12}[4g(\theta (\nabla_{\theta Z}\theta)(W), \nabla_{\theta Y} \xi)+ 6g(\theta(\nabla_{\theta Z} \theta)(\theta Y) , \nabla_{W} \xi)\\
&- 6g(\theta (\nabla_{\theta W}\theta)(\theta Y) , \nabla_{Z} \xi) +g(R(Z,\theta W) \theta Y , \xi) -g(R(W, \theta Z) \theta Y , \xi)]\\
&-\tfrac{17}{6} g(R(W,Z)Y,\xi)+ \tfrac{13}{12} g(R(\theta W,\theta Z)Y,\xi).
\end{align*}
In this equation, we replace $W$ and $Z$ by $\theta W$ and $\theta Z$, to get
\begin{align*}
&2g((\nabla_{\theta W} \theta) (\theta Z) , (\nabla_Y \theta) (\xi)) + 2g((\nabla_Y \theta)(\theta Z) , (\nabla_{\theta W} \theta) (\xi))\\
&=g(\theta \nabla_{\theta W} \xi ,(\nabla_{\theta Z} \theta) (Y)) +g((\nabla_{\theta W} \theta) (Y) , \theta \nabla_{\theta Z} \xi) \\
&+\tfrac{1}{12}[4g(\theta (\nabla_{\theta^2 Z} \theta) (\theta W), \nabla_{\theta Y} \xi)
+ 6g(\theta (\nabla_{\theta^2 Z} \theta) (\theta Y) , \nabla_{\theta W} \xi)\\
&- 6g(\theta (\nabla_{\theta^2 W} \theta) (\theta Y) , \nabla_{\theta Z} \xi) +g(R(\theta Z,\theta^2 W) \theta Y , \xi) 
+g(R(\theta^2 Z , \theta W)\theta Y , \xi)]\\
&-\tfrac{17}{6} g(R(\theta W,\theta Z)Y,\xi)+ \tfrac{13}{12} g(R(\theta^2 W,\theta^2 Z)Y,\xi),
\end{align*}
summing these two equations yields
\begin{align*}
&2g((\nabla_W \theta) (Z) , (\nabla_Y \theta) (\xi)) + 2g((\nabla_Y \theta) (Z) , (\nabla_W \theta) (\xi))\\
&2g((\nabla_{\theta W} \theta) (\theta Z) , (\nabla_Y \theta) (\xi)) + 2g((\nabla_Y \theta) (\theta Z) , (\nabla_{\theta W} \theta) (\xi))\\
&=g(\theta \nabla_W \xi ,(\nabla_Z \theta) (Y)) +g((\nabla_W \theta) (Y) , \theta \nabla_Z \xi) \\
&+g(\theta \nabla_{\theta W} \xi ,(\nabla_{\theta Z} \theta) (Y)) +g((\nabla_{\theta W} \theta) (Y) , \theta \nabla_{\theta Z} \xi) \\
&+\tfrac{1}{12}[4g(\theta (\nabla_{\theta Z} \theta) (W), \nabla_{\theta Y} \xi)+ 6g(\theta (\nabla_{\theta Z} \theta) (\theta Y) , \nabla_{W} \xi)\\
&+ 4g(\theta (\nabla_{\theta^2 Z} \theta) (\theta W), \nabla_{\theta Y} \xi)+ 6g(\theta (\nabla_{\theta^2 Z} \theta) (\theta Y) , \nabla_{\theta W} \xi)\\
&- 6g(\theta (\nabla_{\theta W} \theta) (\theta Y) , \nabla_{Z} \xi) +g(R(Z,\theta W) \theta Y , \xi) -g(R(W, \theta Z) \theta Y , \xi)\\
&- 6g(\theta (\nabla_{\theta^2 W} \theta) (\theta Y) , \nabla_{\theta Z} \xi) +g(R(\theta Z,\theta^2 W) \theta Y , \xi) 
+g(R(\theta^2 Z , \theta W)\theta Y , \xi)]\\
&-\tfrac{17}{6} g(R(W,Z)Y,\xi)+ \tfrac{13}{12} g(R(\theta W,\theta Z)Y,\xi)\\
&-\tfrac{17}{6} g(R(\theta W,\theta Z)Y,\xi)+ \tfrac{13}{12} g(R(\theta^2 W,\theta^2 Z)Y,\xi),
\end{align*}
which simplifies to
\begin{align*}
&g((\nabla_Z \theta) (Y) , \theta (\nabla_W \xi)) + g((\nabla_Y \theta) (W) , \theta (\nabla_Z \xi))
= -\tfrac{7}{4} [g(R(W,Z)Y,\xi)+ g(R(\theta W,\theta Z)Y,\xi)],
\end{align*}
since
\begin{align*}
g(\theta (\nabla_{\theta^2 Z} \theta) (\theta W), \nabla_{\theta Y} \xi)&= -g(\theta (\nabla_{\theta Z} \theta) ( W), \nabla_{\theta Y} \xi).
\end{align*}
Now use Equation~\eqref{eqp45} to obtain
\begin{align*}
&3 g(R(W,Z)Y,\xi) = 2g( \nabla_{Y} \xi ,\theta (\nabla_{Z} \theta) (W)) + g(\nabla_{ W} \xi ,\theta (\nabla_{Z} \theta) (Y)) \\
&-\tfrac{1}{2} g(R(\theta Z, \theta Y)W,\xi) + g(R(\theta W,\theta Z)Y,\xi) - g(\nabla_{Z} \xi ,\theta (\nabla_{W} \theta) (Y)) 
+\tfrac{1}{2} g(R(\theta W, \theta Y)Z,\xi) ,
\end{align*}
but
\begin{align*}
&\tfrac{1}{2} g(R(\theta W, \theta Y)Z,\xi) \\
&= -\tfrac{1}{2} g(R(W,Y)Z,\xi) -\tfrac{2}{7}[ g((\nabla_Y \theta)(Z) , \theta (\nabla_W \xi)) 
+ g((\nabla_Z \theta) (W) , \theta (\nabla_Y \xi))],
\end{align*}
therefore
\begin{align*}
3 g(R(W,Z)Y,\xi) &= 2g( \nabla_{Y} \xi ,\theta (\nabla_{Z} \theta) (W)) + g(\nabla_{ W} \xi ,\theta (\nabla_{Z} \theta) (Y))
- g(\nabla_{Z} \xi ,\theta (\nabla_{W} \theta) (Y))\\
&\tfrac{1}{2} g(R(Z,Y)W,\xi) +\tfrac{2}{7}[ g((\nabla_Y \theta)(W) , \theta (\nabla_Z \xi)) 
+ g((\nabla_W \theta)(Z) , \theta (\nabla_Y \xi))]\\
&-g(R(W,Z)Y,\xi) -\tfrac{4}{7}[ g((\nabla_Z \theta) (Y) , \theta (\nabla_W \xi)) 
+ g((\nabla_Y \theta) (W) , \theta (\nabla_Z \xi))]\\
&-\tfrac{1}{2} g(R(W,Y)Z,\xi) -\tfrac{2}{7}[ g((\nabla_Y \theta)(Z) , \theta (\nabla_W \xi)) 
+ g((\nabla_Z \theta)(W) , \theta (\nabla_Y \xi))],
\end{align*}
and since
$$ \tfrac{1}{2} g(R(Z,Y)W,\xi) -\tfrac{1}{2} g(R(W,Y)Z,\xi) = -\tfrac{1}{2} g(R(W,Z)Y,\xi),$$
we have
\begin{align*}
3 g(R(W,Z)Y,\xi) &= 2g( \nabla_{Y} \xi ,\theta (\nabla_{Z} \theta)(W)) + g(\nabla_{ W} \xi ,\theta (\nabla_{Z} \theta)(Y))
- g(\nabla_{Z} \xi ,\theta (\nabla_{W} \theta) (Y))\\
& +\tfrac{2}{7}[ g((\nabla_Y \theta)(W) , \theta (\nabla_Z \xi)) 
+ g((\nabla_W \theta)(Z) , \theta (\nabla_Y \xi))]\\
&-\tfrac{4}{7}[ g((\nabla_Z \theta)(Y) , \theta (\nabla_W \xi)) 
+ g((\nabla_Y \theta)(W) , \theta (\nabla_Z \xi))]\\
&-\tfrac{2}{7}[ g((\nabla_Y \theta)(Z) , \theta (\nabla_W \xi)) 
+ g((\nabla_Z \theta)(W) , \theta (\nabla_Y \xi))]\\
&-\tfrac{3}{2} g(R(W,Z)Y,\xi) 
\end{align*}
and this finally simplifies into the expression we need:
\begin{align*}
\tfrac{3}{2} g(R(W,Z)Y,\xi) &= \tfrac{18}{7}g( \nabla_{Y} \xi ,\theta (\nabla_{Z} \theta)(W)) + \tfrac{9}{7} g(\nabla_{W} \xi ,\theta (\nabla_{Z} \theta)(Y)) \\
&-\tfrac{9}{7} g(\nabla_{Z} \xi ,\theta(\nabla_{W}\theta)(Y)).
\end{align*}
We now exploit this formula to show $g ( \nabla_{E_i} \xi , R(E_i, X) \xi) =0$.
First, take $W=F_i$, $Z=X$ and $Y=\nabla_{F_i}\xi$, then 
\begin{align*}
g(\nabla_{\nabla_{F_i}\xi} \xi ,\theta (\nabla_{X} \theta) (F_{i}))&= g( \nabla_{\nabla_{\theta F_i}\xi} \xi ,\theta (\nabla_{X} \theta)(\theta F_{i}))\\
&= - g( \nabla_{\nabla_{F_i}\xi} \xi ,\theta (\nabla_{X} \theta) (F_{i}))\\
&=0 ;\\
g(\nabla_{F_i} \xi ,\theta (\nabla_{X} \theta) (\nabla_{F_i}\xi)) &= g(\nabla_{\theta F_i} \xi ,\theta (\nabla_{X} \theta) (\nabla_{\theta F_i}\xi)) \\
&= g(-\theta \nabla_{F_i} \xi ,\theta (\nabla_{X} \theta) (-\theta \nabla_{F_i}\xi)) \\
&= - g(\nabla_{F_i} \xi ,\theta (\nabla_{X}\theta)(\nabla_{F_i}\xi)) \\
&=0;\\
g(\nabla_{X} \xi ,\theta (\nabla_{F_i} \theta) (\nabla_{F_i}\xi))&= 0,
\end{align*}
since $\tr \nabla \theta \circ \nabla \xi =0$. This proves the proposition for $X\in \mathcal F$. 

For $X=\xi$, we go back to Equation~\eqref{frakB} to have
\begin{align*}
g(R(F_i ,\xi)\xi , \nabla_{F_i}\xi) &= - g(R(\xi, F_i)\xi , \nabla_{F_i}\xi)\\
&= g(\nabla_{F_i}\xi , \nabla_{\nabla_{F_i}\xi}\xi)\\
&= g(\nabla_{\theta F_i}\xi , \nabla_{\nabla_{\theta F_i}\xi}\xi)\\
&= g(-\theta \nabla_{F_i}\xi , \nabla_{-\theta \nabla_{F_i}\xi}\xi)\\
&= g(-\theta \nabla_{F_i}\xi , \theta \nabla_{\nabla_{F_i}\xi}\xi)\\
&= -g(\nabla_{F_i}\xi , \nabla_{\nabla_{F_i}\xi}\xi)\\
&= 0 .
\end{align*}

\end{proof}

\end{document}